\documentclass[a4paper,oneside,11pt]{amsart}
\usepackage{enumerate,stackrel}

\usepackage[english]{babel}
\usepackage{fancyhdr}
\usepackage{amsmath}
\usepackage{mathrsfs}
\usepackage{color}

\usepackage{comment}
\let\oldmarginpar\marginpar
\renewcommand\marginpar[1]{\-\oldmarginpar[\raggedleft\small\sf
#1] {\raggedright\small\sf #1}}
\def\esssup{\operatornamewithlimits{\ess\,\sup}}

\usepackage{amssymb}

\newtheorem{definition}{Definition}[section]

\newtheorem{theorem}[definition]{Theorem}

\newtheorem{corollary}[definition]{Corollary}
\newtheorem{lemma}[definition]{Lemma}
\newtheorem{remark}[definition]{Remark}

\numberwithin{equation}{section}

\def\bkR{{\rm I\kern-.17em R}}

\def\bkN{{\rm I\kern-.17em N}}

\def\bkC{{\rm \kern.24em \vrule width.05em height1.4ex depth-.05ex
     \kern-.26em C}}

\def\bkK{{\rm I\kern-.22em K}}
\def\bkP{{\rm I\kern-.22em P}}

\newcommand{\mcal}{\mathcal}

\newcommand{\bpr}{\begin{proof}}
\newcommand{\epr}{\hfill $\square$ \end{proof}}
\newcommand{\bprof}{\begin{proofof}}
\newcommand{\eprof}{\hfill $\square$ \end{proofof}}
\newcommand{\bproff}{\begin{proofoff}}
\newcommand{\eproff}{\hfill $\square$ \end{proofoff}}
\def\esssup{\operatornamewithlimits{ess\,sup}}

\begin{document}


\title{Duals of limiting interpolation spaces} 



\author{Manvi Grover}
\author{Bohum\'{\i}r Opic}


\address{Manvi Grover, 
Department of Mathematical Analysis, 
Faculty of Mathematics and Physics, Charles University,
 Sokolovsk\'a 83, 186 75 Prague 8, Czech Republic}
\email{grover.manvi94@gmail.com}

\address{Bohum\'{\i}r Opic, 
Department of Mathematical Analysis, 
Faculty of Mathematics and Physics, Charles University,
 Sokolovsk\'a 83, 186 75 Prague 8, Czech Republic}
\email{opic@karlin.mff.cuni.cz}




\begin{abstract}

The aim of the paper is to establish duals of the limiting real interpolation $K$- and $J$-spaces   
$(X_0,X_1)_{0,q,v;K}$ and $(X_0,X_1)_{0,q,v;J}$, where  $(X_0,X_1)$ is a compatible couple of Banach spaces, 
$1\le q<\infty$, $v$ is a slowly varying function on the interval $(0, \infty)$, and the symbols $K$ and $J$ stand for the 
Peetre $K$- and $J$-functionals. In the case of the classical real interpolation method  $(X_0,X_1)_{\theta,q}$, where 
$\theta \in (0, 1)$ and $1\le q < \infty$, this problem was solved 
by Lions and Peetre.


\end{abstract}
 
\thanks{The research has been supported by the grant no. 23-04720S of the Czech Science Foundation and by the Grant Schemes at CU, reg. no. CZ.02.2.69/0.0/0.0/19-073/0016935.} 
\maketitle


\section{Introduction}


Dual spaces play a very important role in numerous 
parts of mathematics. The same can be said about spaces which appear in the real interpolation (cf., e.g., monographs  \cite{BB}, \cite{BL}, \cite{Tri78}, \cite{KPS78}, \cite{BS:IO}, \cite{BK}, \cite{L18}).

In this paper our goal is to describe duals of the limiting real interpolation $K$- and $J$-spaces   
$(X_0,X_1)_{0,q,v;K}$ and $(X_0,X_1)_{0,q,v;J}$, where  $(X_0,X_1)$ is a compatible couple of Banach spaces, 
$1\le q<\infty$, $v$ is a slowly varying function on the interval $(0, \infty)$, and the symbols $K$ and $J$ stand for the 
Peetre $K$- and $J$-functionals. Note that we express these duals both in terms of $K$-spaces and $J$-spaces.

In the case $\theta \in (0, 1)$, $1\le q < \infty$, and when the slowly varying function $v$ satisfies $v(t)=1$ for all $t \in (0, \infty)$, we obtain 
the classical real interpolation spaces  $(X_0,X_1)_{\theta,q}$, which are independent on the fact whether the $K$- or $J$-functional
is used to define them (see, e.g.,  \cite [Chapter 5, Theorem 2.8 (Equivalence theorem)]{BS:IO}). Duals of these spaces have been described by Lions and Peetre in their fundamental paper \cite{LP64} (see also  \cite{Li61}). It has been proved there that $(X_0,X_1)'_{\theta,q}=(X_0',X_1')_{\theta,q\,'}$ (with equivalent norms), where $1/q+1/q\,'=1$.

However, some problems in mathematical analysis have motivated the investigation of the real interpolation with the limiting values $\theta=0$ or $\theta=1$.  
Nowadays, there is a~lot of works, where the limiting real interpolation is studied or applied (see, e.g., \cite{GM86}, \cite{Mil94}, \cite{Dok91}, \cite{EO00:RILFR}, \cite{EOP02:RILF}, \cite{GOT2002Lrrinwsvf}, \cite{CFCKU}, \cite{CFCM10}, \cite{AEEK11}, \cite{FMS12}, \cite{EO14}, \cite{FMS14}, 
\cite{COBOS201543}, \cite{CS}, \cite{CFCM15}, \cite{FMS15}, \cite{CD3}, \cite{CDT3}, \cite{Dok18}, \cite{FFGKR18}, \cite{CFC19}, \cite{AFH20}, \cite{BCFC20}, \cite{NO20}, \cite{OG}, \cite{ODST2019MAMS}, \cite{CFCG24}, \cite{O24}, and the references given there).

In order that the spaces $(X_0,X_1)_{\theta,q}$ are meaningful also with limiting values $\theta=0$ or $\theta=1$ in a general case, one has to extend the given interpolation functor by a convenient weight $v$. This weight $v$ usually belongs 
to the class $SV(0,\infty)$ 
of slowly varying functions on the interval $(0,\infty)$. If $\theta \in (0,1)$, then the corresponding space $(X_0,X_1)_{\theta,q,v}$ is a particular case of an interpolation space 
with a function parameter and, by \cite[Theorem 2.2]{Gu78}, this space is again independent on the fact whether the $K$- or $J$-functional
is used to define it.
 However, if $\theta=0$ or $\theta=1$, then the corresponding space depends on the fact whether the $K$- or $J$- functional is used;  accordingly, the resulting space is denoted by  
$(X_0,X_1)_{\theta,q,v;K}$ or by $(X_0,X_1)_{\theta,q,v;J}$. 

Now a natural question arises: Given  $\theta \in \{0,1\}, q \in [1,\infty]$ and $v \in SV(0,\infty)$, can we describe the space $(X_0,X_1)_{\theta,q,v;K}$ as a $(X_0,X_1)_{\theta,q,w;J}$ space with a convenient $w \in SV(0,\infty)$?

If the weight function $v$ is of logarithmic form, then the answer is given in  \cite{CK} for the case that a pair $(X_0,X_1)$ of Banach spaces $X_0$ and $X_1$  is ordered, and in \cite{CS} and \cite{BCFC20} for a~general pair $(X_0,X_1)$ of Banach spaces $X_0$ and $X_1$. 
In \cite{OG} we have establish conditions under which the limiting $K$-space $(X_0,X_1)_{0,q,b;K}$, involving  
a~slowly varying function $b$, can be described by means of the $J$-space $(X_0,X_1)_{0,q,a;J}$, with a convenient slowly varying function $a$, and we have also solved the reverse problem. It has been shown that if these conditions are not satisfied that the given problem may not have a~solution. Moreover, in \cite{O24} it is assumed that these conditions are not fulfilled and then it is proved that 
$$(X_0,X_1)_{0,q,b;K}= (X_0,X_0+X_1)_{0,q,A;J}=X_0+(X_0,X_1)_{0,q,A;J}$$
and 
$$(X_0,X_1)_{0,q,a;J}= (X_0,X_0\cap X_1)_{0,q,B;K}=X_0\cap (X_0,X_1)_{0,q,B;K},$$
where $A$ and $B$ are convenient weights.  In \cite{O24} also equivalent norms in the mentioned spaces have been determined. 

These resuls play important role in our calculation of duals of spaces $(X_0,X_1)_{\theta,q,v;K}$ and $(X_0,X_1)_{\theta,q,v;J}$ 
with $\theta \in \{0,1\}$.   

Note also that it is sufficient to describe duals of these spaces only if $\theta=0$ since the answer for $\theta=1$ follows from the solution with $\theta=0$ and from the equality   
$(X_0,X_1)_{0,q,v;K}=(X_1,X_0)_{1,q,u;K}$, or $(X_0,X_1)_{0,q,v;J}=(X_1,X_0)_{1,q,u;J}$,  
 where $u(t)=v(1/t)$ for all $t>0$ (which is a consequence of the fact that $K(f,t;X_0,X_1)=tK(f,t^{-1};X_1,X_0)$ if  
$f \in X_0+X_1$ and $t>0$, or 
$J(f,t;X_0,X_1)=tJ(f,t^{-1};X_1,X_0)$ if $f \in X_0\cap X_1$ and $t>0$, and a change of variables).

If $\theta \in (0,1)$, $q \in [1,\infty]$, and $v \in SV(0,\infty)$, then, as was already mentioned above, 
the spaces $(X_0,X_1)_{\theta,q,v;K}$,  
$(X_0,X_1)_{\theta,q,v;J}$ coincides and they are particular cases of interpolation spaces with a function parameter. Consequently, 
if $q \in [1,\infty)$, then duals of these spaces are described in \cite[Theorem 2.4]{P86}. Note that it follows from this theorem that 
if $\theta \in (0,1)$, $q \in [1,\infty)$, and $v \in SV(0,\infty)$, then $(X_0,X_1)'_{0,q,v}=(X_0',X_1')_{0,q\,'\!\!,\tilde{v}}$ 
(with equivalent norms), where $1/q+1/q\,'=1$ and $\tilde{v}(t):=\frac{1}{v(1/t)}$ for all $t>0$.

If $\theta\in \{0,1\}$, $q \in [1,\infty)$, and the weight $v$ is of logarithmic form, then the duals of the spaces $(X_0,X_1)_{\theta,q,v;K}$
have been determined in \cite{CS}. 

The paper is organized as follows: Section 2 contains notation, definitions and preliminaries. In Section~3 we present our main 
results. 
In Section~4 we collect some auxiliary assertions. Sections 5--10 are devoted to proofs of the main results. In Section 11 we mention another method which can be used to prove some of our main assertions. 
\hskip 1cm

\section{Notation, definitions and preliminaries}\label{section2}

For two non-negative expressions ({i.e.} functions or functionals) ${\mcal
A}$, ${\mcal B}$, the symbol ${\mcal A}\lesssim {\mcal B}$ (or
${\mcal A}\gtrsim {\mcal B}$) means that $ {\mcal A}\leq c\, {\mcal
B}$ (or $c\,{\mcal A}\geq {\mcal B}$), where $c$ is a  positive constant independent of
significant quantities involved in $A$ and $B$. If ${\mcal A}\lesssim {\mcal
B}$ and ${\mcal A}\gtrsim{\mcal B}$, we write ${\mcal A}\approx
{\mcal B}$ and say that ${\mcal A}$ and ${\mcal B}$ are equivalent.
Throughout the paper we use the abbreviation $\text{LHS}(*)$
($\text{RHS}(*)$) for the left- (right-) hand side of relation
$(*)$. We adopt the convention that $a/(\infty)=0$, 
$a/0=\infty$ and $(\infty)^a=\infty$ for all $a\in(0,\infty)$. If $q\in[1,\infty]$, the conjugate
number $q\,'$ is given by $1/q+1/q\,'=1$. 
In the whole paper $\|.\|_{q;(c,d)},\, q\in
[1,\infty]$, denotes the usual $L_q$-norm on the
interval $(c,d)\subseteq \bkR$.

The symbol ${\mcal M}^{+}(0,\infty)$ stands for the class of all (Lebesgue-) measurable  functions on the interval $(0, \infty)$, which are non-negative almost everywhere on $(0, \infty)$. We will admit only positive, finite weights, and thus we put 
$$
{\mcal W}(0,\infty):=\{w\in{\mcal M}^{+}(0,\infty): w<\infty \ \mbox{a.e. on}\  (0, \infty)\}.
$$
The symbol $AC(0,\infty)$ is used to denote the set of all absolutely continuous functions on the interval $(0,\infty)$.

If $f$ is a monotone function on the interval $(0,\infty)$, then we put 
$$
f(0):=\lim_{t\rightarrow 0_+} f(t) \qquad {\text and} \qquad
f(\infty):=\lim_{t\rightarrow \infty} f(t).
$$

We say that a positive, finite and
Lebesgue-measurable function $b$ is {\it slowly varying} on
$(0,\infty)$, and write $b\in SV(0,\infty)$, if, for each
$\varepsilon>0$, $t^{\varepsilon}b(t)$ is equivalent to a non-decreasing
function on $(0,\infty)$ and  $t^{-\varepsilon}b(t)$ is equivalent to
a~non-increasing function on $(0,\infty)$. Here we follow the definition of $SV(0,\infty)$
given in \cite{GOT2002Lrrinwsvf}; for other definitions see, for example,
\cite{BGT87:RV}, \cite{EEv04:HOFSE}, \cite{EKP:OISRIQ}, and \cite{Nev02:LKSBRPE}. The
family $SV(0,\infty)$ includes not only powers of iterated logarithms
and the broken logarithmic functions of \cite{EO00:RILFR} but also
such functions as $t\mapsto\exp\left(  \left\vert \log
t\right\vert ^{a}\right)  ,$ $a\in(0,1).$ (The last mentioned
function has the interesting property that it tends to infinity more
quickly than any positive power of the logarithmic function).

\smallskip

We mention some properties of slowly
varying functions. 
\begin{lemma} {\rm (\cite [Lemma 2.1]{OG})} \label{l2.6+}
Let $b,b_1,b_2 \in SV(0,\infty)$ and let $d(t):=b(1/t)$ for all $t>0 .$

\noindent
{\rm {\bf (i)}} Then $\, b_1 b_2, d \in SV(0,\infty)$ and $b^r \in SV(0,\infty)$ 
for each $r \in {\mathbb R}.
$ 

\noindent
{\rm {\bf (ii)}}  If $\,\varepsilon$ and $\kappa$ are positive numbers,
then
there are positive constants $c_\varepsilon$ and $C_\varepsilon $ \\
\hspace*{.8cm}such that
\[
c_\varepsilon \min \{ \kappa^{-\varepsilon}, \kappa^\varepsilon \}
b(t) \le b(\kappa t) \le C_\varepsilon \max \{ \kappa^\varepsilon,
\kappa^{-\varepsilon} \} b(t) \quad \mbox{for every }\, t>0 .
\]

\noindent
{\rm {\bf (iii)}} If $\,\alpha>0$ and $q \in (0,\infty]$, then, for all $t>0,$
\[
\|\tau^{\alpha-1/q} b(\tau)\|_{q,(0,t)} \approx t^\alpha b(t) \quad \mbox{and } \quad 
\|\tau^{-\alpha-1/q} b(\tau)\|_{q,(t,\infty)} \approx t^{-\alpha} b(t). 
\] 

\noindent
{\rm {\bf (iv)}} If $\,q \in (0,\infty]$, and     
$$  
B_0(t):=\|\tau^{-1/q} b(\tau)\|_{q,(0,t)}, \quad B_\infty(t):=\|\tau^{-1/q} b(\tau)\|_{q,(t,\infty)}, \ \ t>0,
$$
then
\begin{equation}\label{400}
b(t) \lesssim B_0(t),\quad b(t) \lesssim B_\infty(t), \ \ \text{for all \ } t>0. 
\end{equation}
Moreover, if $B_i(1) <\infty$, then $B_i$ 
belongs to $SV(0,\infty)$, $i=0, \infty$.

\noindent
{\rm {\bf (v)}} If $\,q \in (0,\infty)$, then
\begin{equation}\label{414}
\limsup_{t\rightarrow 0_+}\dfrac{\|\tau^{-1/q} b(\tau)\|_{q,(0,t)}}{b(t)}=\infty,
\quad  \quad
\limsup_{t\rightarrow \infty}\dfrac{\|\tau^{-1/q} b(\tau)\|_{q,(0,t)}}{b(t)}=\infty,
\end{equation}
\begin{equation}\label{415}
\limsup_{t\rightarrow 0_+}\dfrac{\|\tau^{-1/q} b(\tau)\|_{q,(t,\infty)}}{b(t)}=\infty,
\quad  \quad
\limsup_{t\rightarrow \infty}\dfrac{\|\tau^{-1/q} b(\tau)\|_{q,(t,\infty)}}{b(t)}=\infty.
\end{equation}

\end{lemma}

\begin{remark}\label{remark1}\rm{Note that, by Lemma \ref{l2.6+} (iii), the function $b \in SV(0,\infty)$ is equivalent to 
$\overline{b} \in SV(0,\infty) \cap AC(0,\infty)$ given by $\overline{b}(t):= t^{-1}\int_0^t b(s)\,ds,\  t>0$. Consequently, any 
$b \in SV(0,\infty)$ is equvivalent to a continuous function on the interval $(0,\infty).$
} 
\end{remark}

More properties and examples of slowly varying functions can be
found in \cite[Chapter~V, p. 186]{Zyg57:TS}, \cite{BGT87:RV}, 
\cite{VojislavMaric:RVDE00}, \cite{Nev02:LKSBRPE}, 
\cite{GOT2002Lrrinwsvf} and \cite{GNO10:PotAnal}.

\smallskip
Let $X$ and $Y$ be two Banach spaces. We say that $X$ 
{\it coincides} with $Y$ (and write $X=Y$) if $X$ and $Y$ are
equal in the algebraic and topological sense (their norms are equivalent). Moreover, 
we say that $X$ and $Y$ {\it are identical} (and write $X\equiv Y$) provided that $X$ and $Y$ are
equal in the algebraic sense and $\|\cdot \|_X=\|\cdot \|_Y$. 
The symbol $X\hookrightarrow Y$
 means that $X\subset Y$ and the natural embedding of $X$ in $Y$ is continuous. 
The norm of this embedding is denoted by $\|Id\,\|_{X\rightarrow Y}$. By $X'$ we denote the space of all linear and continuous functionals on the space $X$. 

A pair $(X_0,X_1)$ of Banach spaces $X_0$ and $X_1$ is called a {\it compatible couple} if there is a~Hausdorff topological vector space $\mathcal{X}$ in which each of $X_0$ and $X_1$ is continuously embedded.

If  $(X_0,X_1)$ is a compatible couple, then a Banach space $X$ is said to be an {\it intermediate space} 
between $X_0$ and $X_1$ if $X_0\cap X_1 \hookrightarrow X \hookrightarrow X_0+X_1.$

\begin{definition}\label{defin:1:interpolatonspaces}
Let $(X_0,X_1)$ be a compatible couple. 
\begin{trivlist}
\item[\hspace*{0.5cm}{\rm {\bf (i)}}] The Peetre $K$-functional is defined for each $f\in X_0+X_1$ and 
$t>0$ by
$$
K(f,t;X_0,X_1):=\inf\{\|f_0\|_{X_0}+t\|f_1\|_{X_1} : f=f_0+f_1\},
$$
where the infimum extends over all representations $f=f_0+f_1$ of $f$ with $f_0\in X_0$ and 
$f_1\in X_1$. Sometimes, we denote $K(f,t;X_0,X_1)$ simply by $K(f,t)$.

\item[\hspace*{0.5cm}{\rm {\bf (ii)}}] The Peetre $J$-functional is defined for each $f\in X_0\cap X_1$ and 
$t>0$ by
$$
J(f,t;X_0,X_1):= \max \{ \left\|f\right\|_{X_0}, t\left\|f\right\|_{X_1} \}.
$$
Sometimes, we denote $J(f,t;X_0,X_1)$ simply by $J(f,t)$.

\item[\hspace*{0.5cm}{\rm {\bf (iii)}}] For $0\leq \theta\leq 1$, $1\le q\le\infty$, and $v \in 
 {\mcal W}(0,\infty)$, 
we put
\begin{equation}\label{eq:circ_1}
(X_0,X_1)_{\theta,q,v;K}:=\{f\in X_0+X_1: \|f\|_{\theta,q,v;K}<\infty\},
\end{equation}
where
\begin{equation}\label{eq:circ_2}
 \|f\|_{\theta,q,v;K}\equiv\|f\|_{(X_0,X_1)_{\,\theta,q,v;K}}:=\left\|t^{-\theta-1/q}\,v(t)\,K(f,t;X_0,X_1)\right\|_{q,(0,\infty)}.
\end{equation}

\item[\hspace*{0.5cm}{\rm {\bf {(iv)}}}] Let $0\leq \theta\leq 1$, $1\le q\le\infty$, and let $v \in 
 {\mcal W}(0,\infty)$. The space $(X_0,X_1)_{\theta,q,v;J}$ consists of all $f \in X_0+X_1$ for which there is a strongly measurable function $u :(0,\infty)\rightarrow X_0 \cap X_1$ such that 
\begin{equation}\label{300}
f=\int_0^\infty u(s)\, \frac{ds}{s} \quad  {\rm (} {\text convergence\  in \ } X_0+X_1{\rm )}
\end{equation}
and for which the functional
\begin{equation}\label{301}
\|f\|_{\theta,q,v;J}\equiv\|f\|_{(X_0,X_1)_{\,\theta,q,v;J}}:=\inf \left\|t^{-\theta-1/q}\,v(t)\,J(u(t),t;X_0,X_1)\right\|_{q,(0,\infty)}
\end{equation} 
is finite {\rm (}the infimum extends over all representations \eqref{300} of $f${\rm )}.
\end{trivlist}
\end{definition}

We refer to Lemmas \ref{K} and \ref{J}  mentioned below 
for properties of the $K$-functional and the $J$-functional.

\begin{theorem}{\rm (\cite [Theorem 2.3]{OG})} \label{302}
Let $(X_0,X_1)$ be a compatible couple, $0\leq \theta\leq 1$, $1\le q\le\infty$, and let $v \in  {\mcal W}(0,\infty)$.

\noindent
{\rm {\bf A.}} If 
\begin{equation}\label{304}
\left\|t^{-\theta-1/q}\,v(t) \min\{1,t\}\right\|_{q,(0,\infty)}<\infty, 
\end{equation}
then:

{\rm {\bf (i)}} The space $(X_0,X_1)_{\theta,q,v;K}$ is an intermediate space between $X_0$ and $X_1$, that is,
$$
X_0 \cap X_1 \hookrightarrow (X_0,X_1)_{\theta,q,v;K} \hookrightarrow X_0+X_1.
$$

{\rm {\bf (ii)}} 
The space $ (X_0,X_1)_{\theta,q,v;K}$ is a Banach space.

\noindent
{\rm {\bf B.}}
If condition \eqref{304} is not satisfied, then $(X_0,X_1)_{\theta,q,v;K}=\{0\}.$
\end{theorem}


Note that assertion A of Theorem \ref {302} also follows from \cite [Proposition~3.3.1, p. 338]{BK}.

\begin{theorem}{\rm (\cite [Theorem 2.4]{OG})} \label{305}
Let $(X_0,X_1)$ be a compatible couple, $0\leq \theta\leq 1$, $1\le q\le\infty$, and let $v \in  {\mcal W}(0,\infty)$. 

\noindent
{\rm {\bf A.}} If  
\begin{equation}\label{307}
\left\|t^{\theta-1/q\,'}\,\frac{1}{v(t)} \min\left\{1,\frac{1}{t}\right\}\right\|_{q\,'\!,(0,\infty)}<\infty,
\end{equation}
then\,{\rm :}

{\rm {\bf (i)}} The space $(X_0,X_1)_{\theta,q,v;J}$ is an intermediate space between $X_0$ and $X_1$, that is,
$$
X_0 \cap X_1 \hookrightarrow (X_0,X_1)_{\theta,q,v;J} \hookrightarrow X_0+X_1.
$$

{\rm {\bf (ii)}} The space $ (X_0,X_1)_{\theta,q,v;J}$ is a Banach space.

\smallskip
\noindent
\noindent
{\rm {\bf B.}}
If condition \eqref{307} is not satisfied, 
then the functional 
$\|.\|_{\theta,q,v;J}$ vanishes on $X_0\cap X_1$ and thus it is not a norm provided that 
$X_0\cap X_1\ne \{0\}.$
\end{theorem}

Sometimes K-spaces or J-spaces coincide with their modifications, which we now introduce: Let $(X_0,X_1)$ be a compatible couple, $0\leq \theta\leq 1$, $1\le q\le\infty$, and let $v \in  {\mcal W}(0,\infty)$. Assuming that 
$$
{\text either}\ \ (a,b)=(0,1),\qquad\quad {\text or}\ \ (a,b)=(1, \infty),
$$
we put
\begin{equation}\label{eq:circ_1*}
(X_0,X_1)_{\theta,q,v;K;(a,b)}:=\{f\in X_0+X_1: \|f\|_{\theta,q,v;K;(a,b)}<\infty\},
\end{equation}
where
\begin{equation}\label{eq:circ_2*}
 \|f\|_{\theta,q,v;K;(a,b)}\equiv\|f\|_{(X_0,X_1)_{\,\theta,q,v;K;(a,b)}}:=\left\|t^{-\theta-1/q}\,v(t)\,K(f,t;X_0,X_1)\right\|_{q,(a,b)}.
\end{equation}
Similarly, the space $(X_0,X_1)_{\theta,q,v;J;(a,b)}$ consists of all $f \in X_0+X_1$ for which there is a~strongly measurable function $u :(a, b)\rightarrow X_0 \cap X_1$ such that 
\begin{equation}\label{3007}
f=\int_a^b u(s)\, \frac{ds}{s} \quad  {\rm (} {\text convergence\  in \ } X_0+X_1{\rm )}
\end{equation}
and for which the functional
\begin{equation}\label{3017}
\|f\|_{\theta,q,v;J;(a,b)}\equiv\|f\|_{(X_0,X_1)_{\,\theta,q,v;J;(a,b)}}:=\inf \left\|t^{-\theta-1/q}\,v(t)\,J(u(t),t;X_0,X_1)\right\|_{q,(a, b)}
\end{equation} 
is finite {\rm (}the infimum extends over all representations \eqref{3007} of $f${\rm )}.

\smallskip
To calculate duals of interpolation spaces, first note that given a compatible couple $(X_0,X_1)$, then the assumption 
$$
X_0\cap  X_1 \ \text{is dense in }\ X_0 \ \text{and } \ X_1
$$
ensures that the dual couple $(X_0', X_1')$ is compatible as well.

We shall also use the assertions mentioned in the next remark.

\begin{remark}\label{+capD}
{\rm Let $X_0\cap X_1$ is dense in $X_0$ and $X_1$.
 
{\bf (i)} Then \rm{(cf. \cite[p. 53]{BL})}   
\begin{equation}\label{010}
		 K(f',t;X_0',X_1') = \sup_{f\in X_0\cap X_1}\frac{\mid \langle f',f\rangle\mid}{J(f,t^{-1};X_0,X_1)}\quad \text{for all}\ f'\in 
		X_0'+X_1'\ \text{and } t>0,
\end{equation}
and
\begin{equation}\label{2}
    J(f',t;X_0',X_1')=\sup_{f\in X_0+ X_1}\frac{\mid \langle f',f\rangle\mid}{K(f,t^{-1};X_0,X_1)}\quad \text{for all}\ f'\in 
		X_0'\cap X_1'\ \text{and } t>0,
\end{equation}
where $ \langle \cdot,\cdot \rangle$
 denotes the duality 
between $X_0\cap X_1$ 
 and $X_0'+X_1'$ in \eqref{010}, and between 
 $X_0+X_1$ and $X_0'\cap X_1'$ in \eqref{2}.

{\bf (ii)} Since $K(\cdot,1;X_0,X_1)$ is the norm on $X_0 +X_1$ and $J(\cdot,1;X_0,X_1)$ is the norm on $X_0 \cap X_1$ it follows from part (i) 
that
\begin{equation}\label{99}
(X_0\cap X_1)'=X_0'+ X_1'\qquad {\rm{\text and}} \qquad (X_0+X_1)'=X_0'\cap X_1'   
\end{equation}
(which is the result mentioned in \cite [Theorem 2.7.1, p. 32]{BL}).}
\end{remark}

In the following definition we introduce some weighted variant of the small Lebesgue space $\ell_q(\mathbb{Z})$.

\begin{definition}[the space $\lambda_{\theta ,q,w}$]\label{lambda space}
Let $0\leq \theta \leq 1$, $1\leq q\leq \infty$, and let $w\in\mathcal{W}(0,\infty)$. The space $\lambda_{\theta ,q,w}$ consists of all sequences of real numbers $\{\alpha_m\}_{m \in \mathbb{Z}}$ satisfying 
$$
\lVert \{\alpha_m\} \rVert_{\lambda_{\theta,q,w }}=\lVert \{\alpha_m\}_{m\in \mathbb{Z}} \rVert_{\lambda_{\theta,q,w }} <\infty,
$$
where
$$
 \lVert \{\alpha_m\}_{m\in \mathbb{Z}} \rVert_{\lambda_{\theta,q,w }}:= \left\{
    \begin{array}{ll}
       \left(\sum\limits_{m\in\mathbb{Z}}\,(2^{-m\theta}w(2^m)|\alpha_m|)^q\right)^{1/q} & \mbox{if $1\leq q<\infty$, } \\
        \sup\limits_{m\in\mathbb{Z}}\,2^{-m\theta}w(2^m)|\alpha_m| & \mbox{if $q=\infty$.}
    \end{array}
\right.
$$
\end{definition}



The next remark is a consequence of duality results for $\ell_q(\mathbb{Z})$ spaces.

\begin{remark}\label{dual of lamba space}
{\rm Let $0\leq \theta \leq 1$  and $w\in\mathcal{W}(0,\infty).$} 
 {\rm If $1\le q < \infty$, 
then the space $\lambda_{1-\theta ,q\,'\!,w^{-1}}$ is the dual space of $\lambda_{\theta ,q,w}$  via the duality 
$
\sum_{m\in \mathbb{Z}} 2^{-m} \alpha_m \beta_m\,.
$}
\end{remark}

\section{Main Results}\label{sectionmain}

\begin{theorem}[1.\,duality theorem for K-spaces and $\theta=0$]\label{DT0S}
Let $(X_0,X_1)$ be a compatible couple, $1\leq q<\infty$, and let $X_0\cap X_1$ be dense in $X_0$ and $X_1$. If $b\in SV\,(0,\infty)$ satisfies
\begin{equation}\label{DT0A}
\int_x^\infty t^{-1}b^{\,q}(t)\, dt <\infty  \quad\text{for all \ }x>0, \qquad 
\int_0^\infty t^{-1}b^{\,q}(t)\, dt =\infty,
\end{equation}and $a\in$ SV $(0,\infty)$ is defined by
\begin{equation}\label{b=}
a(x):=b^{-{q}/{q\,'}}(x)\int_x^{\infty}t^{-1}b^q(t)\,dt\quad \text{for all }x>0,
\end{equation}
then 
\begin{equation}\label{D2E}
(X_0,X_1)'_{0,q,b;K}=(X_0',X_1')_{0,q\,'\!\!,\tilde{b};J}=(X_0',X_1')_{0,q\,'\!\!,\tilde{a};K}, 
\end{equation} 
where $\tilde{b}(x):=\frac{1}{b({1}/{x})}$ and $\tilde{a}(x):=\frac{1}{a({1}/{x})}$ for all $x>0$.
\end{theorem}

\smallskip
(The proof of this theorem is in Section 5.) 

\medskip
There are the following two counterparts of Theorem \ref{DT0S}.
\begin{theorem}[1. duality theorem for J-spaces and $\theta=0$]
\label{DTJ1}
Let $(X_0, X_1)$ be a compatible couple, $1<q< \infty$, and let $X_0\cap X_1$ be dense in $X_0$ and $X_1$.
 If $a \in SV(0,\infty)$ satisfies
\begin{equation}\label{prop_slow_var_funct_aJ}
\int_0^x t^{-1}a^{-q\,'}(t)\, dt <\infty  \quad\text{for all \ }x>0, \qquad 
\int_0^\infty t^{-1}a^{-q\,'}(t)\, dt =\infty,
\end{equation}
and $b \in  SV(0,\infty)$ is defined by
\begin{equation}\label{def_slow_var_funct_bJ}
b(x):=a^{-q\,'\!/q}(x)\Big(\int_0^x t^{-1}a^{-q\,'}(t)\, dt \Big)^{-1} \quad\text{for all \ }x>0,
\end{equation}
then 
\begin{equation}\label{1J}
(X_0, X_1)'_{0, q, a; J}=(X'_0, X'_1)_{0, q\,'\!\!, \tilde{a}; K}=(X'_0, X'_1)_{0, q\,'\!\!, \tilde{b}; J},
\end{equation}
where $\tilde{a}(x):=\frac{1}{a({1}/{x})}$ and $\tilde{b}(x):=\frac{1}{b({1}/{x})}$ for all $x>0$.
\end{theorem}

(The proof of this theorem is in Section \ref{Pr3Main}.) 

\medskip

\begin{theorem}[2. duality theorem for J-spaces and $\theta=0$]
\label{DTJ11}
Let $(X_0, X_1)$ be a compatible couple and let $X_0\cap X_1$ be dense in $X_0$ and $X_1$. 
Assume that $a \in SV(0,\infty)\cap AC(0,\infty)$ satisfies
\begin{equation}\label{prop_slow_var_funct_a*1}
a \text{\ is strictly decreasing},\quad a(0)=\infty, \quad a(\infty)=0.
\end{equation}
If $b \in  SV(0,\infty)$ and 
\begin{equation}\label{def_slow_var_funct_b*1}
b(x):= -x \,a\, '(x) \quad\text{for a.a. }x>0,
\end{equation}
then
\begin{equation}\label{1J1}
(X_0, X_1)'_{0, 1, a; J}=(X'_0, X'_1)_{0, \infty, \tilde{a}; K}=(X'_0, X'_1)_{0, \infty, \tilde{b}; J},
\end{equation}
where $\tilde{a}(x):=\frac{1}{a({1}/{x})}$ and $\tilde{b}(x):=\frac{1}{b({1}/{x})}$ for all $x>0$.
\end{theorem}

(The proof of this theorem is in Section \ref{Pr4Main}.)

\medskip
There is the following variant of Theorem \ref{DT0S}.

\begin{theorem}[2.\,duality theorem for K-spaces and $\theta=0$]\label{DT0S1}
Let $(X_0,X_1)$ be a compatible couple, $1\leq q<\infty$, and let $X_0\cap X_1$ be dense in $X_0$ and $X_1$. If $b\in SV\,(0,\infty)$ satisfies
\begin{equation}\label{DT0A1}
\int_0^\infty t^{-1}b^{\,q}(t)\, dt <\infty,
\end{equation}
the function $B\in SV\,(0,\infty)$ is given by 
\begin{equation}\label{3131}
B(x):=b(x) \ \ \text{if} \ x\in [1,\infty)\quad \text{and}\quad B(x)\approx \beta(x) \ \ \text{if} \ x\in (0,1),
\end{equation}
where $\beta\in SV\,(0,\infty)$ is such that
\begin{equation}\label{3141}
      \int_x^\infty t^{-1} \beta^q(t)\, dt <\infty \quad \text{for all} \ x>0, \quad \quad   \int_0^\infty t^{-1} \beta^q(t)\, 
			dt=\infty, 
\end{equation}
and 
\begin{equation}\label{B1=}
A(x):=B^{-{q}/{q\,'}}(x)\int_x^{\infty}t^{-1}B^q(t)\,dt\quad \text{for all }x>0,
\end{equation}
then 
\begin{align}\label{D2E1}
(X_0,X_1)'_{0,q,b;K}=&(X_0', X_1')_{0,q\,'\!\!,\tilde{b};J}\\
=&(X_0',X_0'\cap X_1')_{0,q\,'\!\!,\tilde{B};J}=
(X_0',X_0'\cap X_1')_{0,q\,'\!\!,\tilde{A};K}\notag\\
=&X_0'\cap(X_0',X_1')_{0,q\,'\!\!,\tilde{B};J}=
X_0'\cap (X_0',X_1')_{0,q\,'\!\!,\tilde{A};K},\notag 
\end{align} 
where $\tilde{b}(x):=\frac{1}{b({1}/{x})}, \tilde{B}(x):=\frac{1}{B({1}/{x})}$, and 
$\tilde{A}(x):=\frac{1}{A({1}/{x})}$ for all $x>0$.
\end{theorem}


\begin{remark}\label{remark100*}\rm{Under the assumptions of Theorem \ref{DT0S1},
\begin{equation}\label{101*}
(X_0,X_1)_{0,q,b;K}=(X_0,X_1)_{0,q,b;K;(1,\infty)},
\end{equation}
\begin{equation}\label{1002*}
(X_0', X_1')_{0,q\,'\!\!,\tilde{b};J}=(X_0', X_1')_{0,q\,'\!\!,\tilde{b};J;(0,1)},
\end{equation}
\begin{equation}\label{1003*}
(X_0',X_0'\cap X_1')_{0,q\,'\!\!,\tilde{B};J}=(X_0',X_0'\cap X_1')_{0,q\,'\!\!,\tilde{B};J;(0,1)},
\end{equation}
\begin{equation}\label{1004*}
(X_0',X_0'\cap X_1')_{0,q\,'\!\!,\tilde{A};K}=(X_0',X_0'\cap X_1')_{0,q\,'\!\!,\tilde{A};K;(0,1)}.
\end{equation}
}
\end{remark}
(The proofs of Theorem \ref{DT0S1} and Remark \ref{remark100*} are in Section \ref{Pr2Main}.) 

\medskip
There is the following variant of Theorem \ref{DTJ1}.
\begin{theorem}[3. duality theorem for J-spaces and $\theta=0$]
\label{DTJ111}
Let $(X_0, X_1)$ be a compatible couple, $1<q<\infty$, and let $X_0\cap X_1$ be dense in $X_0$ and $X_1$. 
If $a \in SV(0,\infty)$ satisfies
\begin{equation}\label{akon}
\int_0^{\infty} t^{-1}a^{-q\,'}(t)\, dt <\infty,
\end{equation}
the function $A\in SV(0,\infty)$ is given by
\begin{equation}\label{Afce}
A(x):=a(x)\quad\text{if \ } x\in (0,1], \qquad A(x)\approx \alpha(x) \quad\text{if \ } x\in (1, \infty),
\end{equation}
where the function $\alpha \in SV(0,\infty)$ is such that
\begin{equation}\label{alfa}
\int_0^x t^{-1}\alpha^{-q\,'}(t)\, dt <\infty\quad\text{for all \ }x>0, \qquad 
\int_0^\infty t^{-1}\alpha^{-q\,'}(t)\, dt =\infty,
\end{equation}
and
\begin{equation}\label{Bfce}
B(x):=A^{-q\,'\!/q}(x) \left(\int_0^x t^{-1}A^{-q\,'}(t)\,dt\right)^{-1} \quad\text{for all \ }x>0,
\end{equation}
then
\begin{align}\label{dual}
(X_0, X_1)'_{0, q, a; J}=&(X_0',X_1')_{0,q\,'\!\!,\tilde{a};K}\\
=&(X_0',X_0'+X_1')_{0,q\,'\!\!,\tilde{A};K}=(X'_0, X'_0+X'_1)_{0, q\,'\!\!, \tilde{B};J}\notag\\
=&X_0'+(X_0',X_1')_{0,q\,'\!\!,\tilde{A};K}=X'_0+ (X'_0,X'_1)_{0, q\,'\!\!, \tilde{B};J},\notag
\end{align}
where $\tilde{a}(x):=\frac{1}{a({1}/{x})}, \tilde{A}(x):=\frac{1}{A({1}/{x})}$, and $\tilde{B}(x):=\frac{1}{B({1}/{x})}$ for all $x>0$.
\end{theorem}

\begin{remark}\label{remark100J*}\rm{Under the assumptions of Theorem \ref{DTJ111},
\begin{equation}\label{101J*}
(X_0,X_1)_{0,q,a;J}=(X_0,X_1)_{0,q,a;J;(0,1)},
\end{equation}
\begin{equation}\label{1002J*}
(X_0',X_1')_{0,q\,'\!\!,\tilde{a};K}=(X_0',X_1')_{0,q\,'\!\!,\tilde{a};K;(1,\infty)},
\end{equation}
\begin{equation}\label{1003J*}
(X_0',X_0'+X_1')_{0,q\,'\!\!,\tilde{A};K}=(X_0',X_0'+X_1')_{0,q\,'\!\!,\tilde{A};K;(1,\infty)},
\end{equation}
\begin{equation}\label{1004J*}
(X'_0, X'_0+X'_1)_{0, q\,'\!\!, \tilde{B}; J}=(X'_0, X'_0+X'_1)_{0, q\,'\!\!, \tilde{B}; ;(1,\infty)}.
\end{equation}
}
\end{remark}

\smallskip
(The proofs of Theorem \ref{DTJ111} and Remark \ref{remark100J*} are Section \ref{Pr5Main}.)

\medskip
The next assertion is a variant of Theorem \ref{DTJ11}.

\begin{theorem}[4. duality theorem for J-spaces and $\theta=0$]\label{DTJ111.1}
Let $(X_0, X_1)$ be a compatible couple and let $X_0\cap X_1$ be dense in $X_0$ and $X_1$. 
Assume that $a \in SV(0,\infty)\cap AC(0,\infty)$ satisfies
\begin{equation}\label{akon.1}
a \text{  is strictly decreasing}, \qquad a(0)=\infty, \qquad a(\infty)>0,
\end{equation}
the function $A \in SV(0,\infty)$ be given by
\begin{equation}\label{Afce.1}
A(x):=a(x)\quad\text{if \ } x\in (0,1], \qquad A(x):=c\,\alpha(x) \quad\text{if \ } x\in (1,\infty),
\end{equation}
where $\alpha \in SV(0,\infty)\cap AC(0,\infty)$ is such that 
\begin{equation}\label{alfa.1}
\alpha \text{ is strictly decreasing},\qquad  \alpha(0)=\infty,\qquad  \alpha(\infty)=0,
\end{equation}
and $c$ is a positive constant chosen in such a way that $A \in AC(0,\infty)$. If $B \in SV(0,\infty)$ and 
\begin{equation}\label{Bfce.1}
B(x):= -x \,A\, '(x) \quad\text{for a.a. }x>0,
\end{equation}
then 
\begin{align}\label{vnoreni.1}
(X_0, X_1)'_{0, 1, a; J}&=(X_0',X_1')_{0,\infty,\tilde{a};K}\\
&=(X'_0, X'_0+X'_1)_{0, \infty, \tilde{A}; K}=(X'_0, X'_0+X'_1)_{0, \infty, \tilde{B}; J}\notag \\
&=X'_0+(X_0',X_1')_{0,\infty,\tilde{A};K}=X'_0+(X_0',X_1')_{0,\infty,\tilde{B};J},\notag
\end{align}
where $\tilde{a}(x):=\frac{1}{a({1}/{x})}, \tilde{A}(x):=\frac{1}{A({1}/{x})}$, and $\tilde{B}(x):=\frac{1}{B({1}/{x})}$ for all $x>0$.
\end{theorem}

\begin{remark}\label{remark100J*4}\rm{Under the assumptions of Theorem \ref{DTJ111.1},
\begin{equation}\label{101J*4}
(X_0,X_1)_{0,1,a;J}=(X_0,X_1)_{0,1,a;J;(0,1)},
\end{equation}
\begin{equation}\label{1002J*4}
(X_0',X_1')_{0,\infty,\tilde{a};K}=(X_0',X_1')_{0,\infty,\tilde{a};K;(1,\infty)},
\end{equation}
\begin{equation}\label{1003J*4}
(X_0',X_0'+X_1')_{0,\infty,\tilde{A};K}=(X_0',X_0'+X_1')_{0,\infty,\tilde{A};K;(1,\infty)}.
\end{equation}
Moreover, if $\|\tilde{B}(t)\|_{\infty, (1,e)}<\infty$, then also
\begin{equation}\label{1004J*4}
(X'_0, X'_0+X'_1)_{0,\infty,\tilde{B}; J}=(X'_0, X'_0+X'_1)_{0, \infty, \tilde{B}; J;(1,\infty)}.
\end{equation}
}
\end{remark}
 
\smallskip
(The proofs of Theorem \ref{DTJ111.1} and Remark \ref{remark100J*4} are in Section \ref{Pr6Main}.)

\begin{remark}\label{remark2}\rm{In Theorem \ref{DTJ11} we assume that $a \in SV(0,\infty) \cap AC(0,\infty)$. By 
Remark~\ref{remark1}, any function $a \in SV(0,\infty)$ is equvivalent to the function 
$\overline{a} \in SV(0,\infty) \cap AC(0,\infty)$ given by $\overline{a}(x):= x^{-1}\int_0^x a(t)\,dt,\  x>0$. Moreover, one can 
prove that if the function $a$ satisfies \eqref{prop_slow_var_funct_a*1}, then \eqref{prop_slow_var_funct_a*1} also holds with $a$ replaced by $\overline{a}$. Thus, Theorem \ref{DTJ11} remains true if the assumption $a \in SV(0,\infty) \cap AC(0,\infty)$ is replaced 
only by $a \in SV(0,\infty)$ provided that in \eqref{def_slow_var_funct_b*1} we write $\overline{a}\,'$ instead of $a'$.

Similar remark can be made about Theorem \ref{DTJ111.1}.} 
\end{remark}

\medskip

\section{Auxiliary assertions}

In the next  two lemmas some basic properties of the $K$- and $J$- functionals are summarized.

\begin{lemma} {\rm (cf. \cite [Proposition 1.2, p. 294]{BS:IO})} \label{about K}
If $(X_0,X_1)$ is a compatible couple, then, for each $f \in X_0+X_1$, the $K$-functional $K(f,t;X_0,X_1)$ is a nonnegative concave function of $t>0$, and 
\begin{equation}\label{K}
K(f,t;X_0,X_1)=t\,K(f,t^{-1};X_1,X_0) \quad \text{for all \ } t>0.
\end{equation}
In particular, $K(f,t;X_0,X_1)$ is non-decreasing on $(0,\infty)$ and $K(f,t;X_0,X_1)/t$ is non-increasing 
on $(0,\infty).$ 
\end{lemma}

\begin{lemma} {\rm (cf. \cite [Lemma 3.2.1, p. 42]{BL})}\label{about J}
If $(X_0,X_1)$ is a compatible couple, then, for each $f \in X_0\cap X_1$, the $J$-~functional $J(f,t;X_0,X_1)$ is a nonnegative convex function of $t>0$, and 
\begin{equation}\label{J}
J(f,t;X_0,X_1)=t\,J(f,t^{-1};X_1,X_0)\quad \text{for all \ } t>0,
\end{equation}
\begin{equation}\label{JJ}
J(f,t;X_0,X_1)\le \max\{1, t/s\}\,J(f,s;X_0,X_1) \quad \text{for all \ } t, s>0,
\end{equation}
\begin{equation}\label{KJ}
K(f,t;X_0,X_1)\le \min\{1, t/s\}\,J(f,s;X_0,X_1)\quad \text{for all \ } t, s>0.
\end{equation}
In particular, $J(f,t;X_0,X_1)$ is non-decreasing on $(0,\infty)$ and $J(f,t;X_0,X_1)/t$ is non-increasing 
on $(0,\infty).$ 
\end{lemma}

The following result describes when the $K$-space $(X_0, X_1)_{\theta, q, b; K}$ with the limiting value 
$\theta=0$ and $b \in SV(0,\infty)$ coincides with the $J$-space $(X_0, X_1)_{\theta, q, a; J}$ with a  convenient $a \in SV(0,\infty).$

\begin{theorem}{\rm (cf. \cite [Theorem 3.1]{OG})} 
\label{equivalence theorem1}
Let $(X_0, X_1)$ be a compatible couple and $1\le q <\infty$. If $b \in SV(0,\infty)$ satisfies
\begin{equation}\label{prop_slow_var_funct_b1}
\int_x^\infty t^{-1}b^{\,q}(t)\, dt <\infty  \quad\text{for all \ }x>0, \qquad 
\int_0^\infty t^{-1}b^{\,q}(t)\, dt =\infty,
\end{equation}
and $a \in  SV(0,\infty)$ is such that 
\begin{equation}\label{def_slow_var_funct_a1}
a(x)\approx b^{-q/q\,'}(x)\int_x^\infty t^{-1}b^{\,q}(t)\, dt  \quad\text{for a.a.\ }x>0,
\end{equation}
then 
\begin{equation}\label{K_space=J_space1}
(X_0, X_1)_{0, q, b; K}=(X_0, X_1)_{0, q, a; J}.
\end{equation}
\end{theorem}

\smallskip
\begin{theorem}{\rm (cf. \cite [Theorem 3.2]{OG})} \label{density theorem1}
Let $(X_0, X_1)$ be a compatible couple and $1\le q<\infty$. If $b \in SV(0,\infty)$ satisfies
\eqref{prop_slow_var_funct_b1}, 
then  the space $X_0\cap X_1$ is dense in $(X_0, X_1)_{0, q, b; K}$.
\end{theorem}

\smallskip
There are the following counterparts of the previous two assertions.

\begin{theorem}{\rm (cf. \cite [Theorem 3.3]{OG})} 
\label{equivalence theorem}
Let $(X_0, X_1)$ be a compatible couple and $1<q\le \infty$. If $a \in SV(0,\infty)$ satisfies
\begin{equation}\label{prop_slow_var_funct_a}
\int_0^x t^{-1}a^{-q\,'}(t)\, dt <\infty  \quad\text{for all \ }x>0, \qquad 
\int_0^\infty t^{-1}a^{-q\,'}(t)\, dt =\infty,
\end{equation}
and $b \in  SV(0,\infty)$ is such that 
\begin{equation}\label{def_slow_var_funct_b}
b(x) \approx a^{-q\,'\!/q}(x)\Big(\int_0^x t^{-1}a^{-q\,'}(t)\, dt \Big)^{-1} \quad\text{for a.a.\ }x>0,
\end{equation}
then 
\begin{equation}\label{K_space=J_space}
(X_0, X_1)_{0, q, a; J}=(X_0, X_1)_{0, q, b; K}.
\end{equation}
\end{theorem}

\smallskip
\begin{theorem}{\rm (cf. \cite [Theorem 3.4]{OG})} \label{density theorem}
Let $(X_0, X_1)$ be a compatible couple and $1<q<\infty$. If \,$a \in SV(0,\infty)$ satisfies
\eqref{prop_slow_var_funct_a}, 
then  the space $X_0\cap X_1$ is dense in $(X_0, X_1)_{0, q, a; J}$.
\end{theorem}

\smallskip
The next two assertions are complements of Theorems \ref{equivalence theorem1} and \ref{equivalence theorem}.

\begin{theorem}{\rm (\cite [Theorem 9.1]{OG})}
\label{equivalence theorem1*}
Let $(X_0, X_1)$ be a compatible couple. If $b \in SV(0,\infty)\cap AC(0,\infty)$ satisfies
\begin{equation}\label{prop_slow_var_funct_b1*}
b \text{\ is strictly decreasing},\quad b(0)=\infty, \quad b(\infty)=0,
\end{equation}
and if 
\begin{equation}\label{def_slow_var_funct_a1*}
a(x):= \frac{b^{\,2}(x)}{x (-b\, '(x))} \quad\text{for a.a. }x>0,
\end{equation}
then 
\begin{equation}\label{K_space=J_space1*}
(X_0, X_1)_{0, \infty, b; K}=(X_0, X_1)_{0, \infty, a; J}.
\end{equation}
\end{theorem}

\smallskip
\begin{theorem}{\rm(\cite [Theorem 10.1]{OG})}
\label{equivalence theorem1**}
Let $(X_0, X_1)$ be a compatible couple. If $a \in SV(0,\infty)\cap AC(0,\infty)$ satisfies
\begin{equation}\label{prop_slow_var_funct_a*}
a \text{\ is strictly decreasing},\quad a(0)=\infty, \quad a(\infty)=0,
\end{equation}
and if
\begin{equation}\label{def_slow_var_funct_b*}
b(x):= -x \,a\, '(x) \quad\text{for a.a. }x>0,
\end{equation}
then 
\begin{equation}\label{K_space=J_space1**}
(X_0, X_1)_{0, 1, a; J}=(X_0, X_1)_{0, 1, b; K}.
\end{equation}
\end{theorem}

\smallskip
\begin{theorem}{\rm(\cite [Theorem 10.3]{OG})}\label{density theorem*}
Let $(X_0, X_1)$ be a compatible couple. If the function $a$ 
satisfies the assumptions of {\rm  Theorem \ref{equivalence theorem1**}}, then 
$X_0\cap X_1$ is dense in $(X_0, X_1)_{0, 1
, a; J}$.
\end{theorem}

\smallskip
Another complement of Theorem \ref{equivalence theorem1} is the next result.

\begin{theorem}{\rm(\cite [Theorem 1.11]{O24})}
\label{KS}
Let $(X_0,X_1)$ be a compatible couple and $1\le q<\infty$. If $b\in SV\,(0,\infty)$ satisfies
\begin{equation}\label{KSB}
\int_0^\infty t^{-1}b^{\,q}(t)\, dt <\infty,  
\end{equation}
the function $B\in SV\,(0,\infty)$ is given by 
\begin{equation}\label{3131*}
B(x):=b(x) \ \ \text{if} \ x\in [1,\infty),\qquad 
B(x)\approx \beta(x) \ \ \text{if} \ x\in (0,1),
\end{equation}
where $\beta\in \, SV (0,\infty)$ is such that
\begin{equation}\label{3141*}
      \int_x^\infty t^{-1} \beta^q(t)\, dt <\infty \quad \text{for all} \ x>0, \quad \quad   
			\int_0^\infty t^{-1} \beta^q(t)\, dt=\infty, 
\end{equation}
and 
\begin{equation}\label{B1=*}
A(x):=B^{-{q}/{q\,'}}(x)\int_x^{\infty}t^{-1}B^q(t)\,dt\quad \text{for all }x>0,
\end{equation}
then 
\begin{equation}\label{D2E1*}
(X_0,X_1)_{0,q,b;K}=(X_0, X_0+X_1)_{0,q,B;K}=(X_0,X_0+X_1)_{0,q,A;J}.
\end{equation} 
\end{theorem}

\begin{remark}{\rm(\cite [Remark 1.13]{O24})}\label{remark100} \rm{Under the assumptions of Theorem \ref{KS},
\begin{equation}\label{101}
(X_0,X_1)_{0,q,b;K}=(X_0,X_1)_{0,q,b;K;(1,\infty)},
\end{equation}
\begin{equation}\label{1002}
(X_0,X_0+X_1)_{0,q,B;K}=(X_0,X_0+X_1)_{0,q,B;K;(1,\infty)},
\end{equation}
\begin{equation}\label{1003}
(X_0,X_0+X_1)_{0,q,A;J}=(X_0,X_0+X_1)_{0,q,A;J;(1,\infty)}.
\end{equation}
}
\end{remark}

\begin{corollary}{\rm(\cite [Corollary 1.14]{O24})}\label{cor2/2}
Under the assumptions of {\rm Theorem \ref{KS}},
\begin{equation}\label{corres2/2}
(X_0, X_1)_{0, q, b; K}=X_0 + (X_0,X_1)_{0, q, B; K}=X_0 + (X_0,X_1)_{0, q, A; J}.
\end{equation}
\end{corollary}

\smallskip
Further complement of Theorem \ref{equivalence theorem} reads as follows.

\begin{theorem}{\rm(\cite [Theorem 1.15]{O24})}\label{JS11}
Let $(X_0, X_1)$ be a compatible couple and $1<q\le\infty$.
If $a \in SV(0,\infty)$ satisfies
\begin{equation}\label{akon111}
\int_0^{\infty} t^{-1}a^{-q\,'}(t)\, dt <\infty,
\end{equation}
the function $A\in SV(0,\infty)$ is given by
\begin{equation}\label{Afce111}
A(x):=a(x) \quad\text{if \ } x\in (0,1], \qquad 
A(x)\approx \alpha(x)\quad\text{if \ } x\in (1, \infty),
\end{equation}
where the function $\alpha \in SV(0,\infty)$ is such that
\begin{equation}\label{alfa111}
\int_0^x t^{-1}\alpha^{-q\,'}(t)\, dt <\infty\quad\text{for all \ }x>0, \qquad 
\int_0^\infty t^{-1}\alpha^{-q\,'}(t)\, dt =\infty,
\end{equation}
and
\begin{equation}\label{Bfce111}
B(x):=A^{-q\,'\!/q}(x) \left(\int_0^x t^{-1}A^{-q\,'}(t)\,dt\right)^{-1} \quad\text{for all \ }x>0,
\end{equation}
then
\begin{equation}\label{dual*}
(X_0, X_1)_{0, q, a; J}=(X_0, X_0\cap X_1)_{0,q,A;J}=(X_0, X_0\cap X_1)_{0, q, B; K}.
\end{equation}
\end{theorem}

\begin{remark}{\rm(\cite [Remark 1.17]{O24})}
\label{remark100J}\rm{Note that, under the assumptions of Theorem \ref{JS11},
\begin{equation}\label{101J}
(X_0,X_1)_{0,q,a;J}=(X_0,X_1)_{0,q,a;J;(0,1)},
\end{equation}
\begin{equation}\label{1002J}
(X_0,X_0\cap X_1)_{0,q,A;J}=(X_0,X_0\cap X_1)_{0,q,A;J;(0,1)},
\end{equation}
\begin{equation}\label{1003J}
(X_0,X_0\cap X_1)_{0,q,B;K}=(X_0,X_0\cap X_1)_{0,q,B;K;(0,1)}.
\end{equation}
}
\end{remark}

\begin{corollary}{\rm(\cite [Corollary 1.18]{O24})}\label{cor2}
Under the assumptions of {\rm Theorem \ref{JS11}},
\begin{equation}\label{corres2}
(X_0, X_1)_{0, q, a; J}=X_0 \cap (X_0,X_1)_{0, q, A; J}=X_0 \cap (X_0,X_1)_{0, q, B; K}.
\end{equation}
\end{corollary}

\smallskip

The following result is a complement of Theorem \ref{equivalence theorem1**}.

\begin{theorem}{\rm (\cite[Theorem 1.22]{O24})}
\label{equivalence theorem1**a}
Let $(X_0, X_1)$ be a compatible couple and let $a \in SV(0,\infty)\cap AC(0,\infty)$ satisfy
\begin{equation}\label{prop_slow_var_funct_a**}
a \text{\ is strictly decreasing},\quad a(0)=\infty, \quad a(\infty)>0.
\end{equation}
Assume that the function $A\in SV(0,\infty)$ is given by 
\begin{equation}\label{3131*A}
A(x):=a(x) \ \ \text{if} \ x\in (0,1],\qquad 
A(x):=c\,\alpha(x) \ \ \text{if} \ x\in (1,\infty),
\end{equation}
where $\alpha\in \, SV(0,\infty)\cap AC(0,\infty)$ is such that
\begin{equation}\label{bq1}
\alpha \text{\ is strictly decreasing},\quad \alpha(0)=\infty, \quad \alpha(\infty)=0,      
\end{equation}
and $c$ is a positive constant chosen in such a way that $A \in AC(0,\infty)$. If 
\begin{equation}\label{difeqB}
B(x):=-x \,A\, '(x) \quad\text{for a.a. }x>0, 
\end{equation}
then
\begin{equation}\label{D2E1*S}
(X_0,X_1)_{0,1,a;J}=(X_0, X_0\cap X_1)_{0,1,A;J}=(X_0,X_0\cap X_1)_{0,1,B;K}.
\end{equation} 
\end{theorem}

\smallskip
\begin{theorem} {\rm (\cite[Theorem 1.23]{O24})}
\label{density theorem**}
Let $(X_0, X_1)$ be a compatible couple. If the function~$a$ 
satisfies the assumptions of {\rm  Theorem \ref{equivalence theorem1**a}}, then 
$X_0\cap X_1$ is dense in $(X_0, X_1)_{0, 1, a; J}$.
\end{theorem}

\smallskip
\begin{remark}{\rm (\cite[Remark 1.24]{O24})}
\label{remark100J1} \rm{Note that, under the assumptions of Theorem \ref{equivalence theorem1**a},
\begin{equation}\label{101J1}
(X_0,X_1)_{0,1,a;J}=(X_0,X_1)_{0,1,a;J;(0,1)},
\end{equation}
\begin{equation}\label{1002J1}
(X_0,X_0\cap X_1)_{0,1,A;J}=(X_0,X_0\cap X_1)_{0,1,A;J;(0,1)},
\end{equation}
\begin{equation}\label{1003J1}
(X_0,X_0\cap X_1)_{0,1,B;K}=(X_0,X_0\cap X_1)_{0,1,B;K;(0,1)}.
\end{equation}
} 
\end{remark}

\smallskip
\begin{corollary}{\rm(\cite [Corollary 1.25]{O24})}\label{cor21} 
Under the assumptions of {\rm Theorem \ref{equivalence theorem1**a}},
\begin{equation}\label{corres21}
(X_0, X_1)_{0, 1, a; J}=X_0 \cap (X_0,X_1)_{0, 1, A; J}=X_0 \cap (X_0,X_1)_{0, 1, B; K}.
\end{equation}
\end{corollary}

\smallskip

There is the following counterpart of Theorem \ref{equivalence theorem1**a}.

\begin{theorem}[{\cite[Theorem 1.19]{O24}}]
\label{*equivalence theorem1*}
Let $(X_0, X_1)$ be a compatible couple and let $b \in SV(0,\infty)\cap AC(0,\infty)$ satisfy
\begin{equation}\label{*prop_slow_var_funct_b1*}
b \text{\ is strictly decreasing},\quad b(0)<\infty, \quad b(\infty)=0.
\end{equation}
Assume that the function $B\in SV(0,\infty)$ is given by 
\begin{equation}\label{3131*2}
B(x):=b(x) \ \ \text{if} \ x\in [1,\infty),\qquad 
B(x):=c\,\beta(x) \ \ \text{if} \ x\in (0,1),
\end{equation}
where $\beta\in \, SV(0,\infty)\cap AC(0,\infty)$ is such that
\begin{equation}\label{aqinfty}
\beta \text{\ is strictly decreasing},\quad \beta(0)=\infty, \quad \beta(\infty)=0,      
\end{equation}
and $c$ is a positive constant chosen in such a way that $B \in AC(0,\infty)$. 
If 
\begin{equation}\label{difeq}
A(x):= \frac{B^{\,2}(x)}{x (-B\, '(x))} \quad\text{for a.a. }x>0,
\end{equation}
then
\begin{equation}\label{*D2E1*2}
(X_0,X_1)_{0,\infty,b;K}=(X_0, X_0+X_1)_{0,\infty,B;K}=(X_0,X_0+X_1)_{0,\infty,A;J}.
\end{equation} 
\end{theorem}

\smallskip
\begin{remark}{\rm (\cite[Remark 1.20]{O24})}
\label{remark100*4}\rm{Note that, under the assumptions of Theorem \ref{*equivalence theorem1*},
\begin{equation}\label{101*4}
(X_0,X_1)_{0,\infty,b;K}=(X_0,X_1)_{0,\infty,b;K;(1,\infty)},
\end{equation}
\begin{equation}\label{1002*4}
(X_0,X_0+X_1)_{0,\infty,B;K}=(X_0,X_0+X_1)_{0,\infty,B;K;(1,\infty)}.
\end{equation}
Moreover, if $\|A(t)\|_{\infty, (1,e)}<\infty$, then also
\begin{equation}\label{1003*4}
(X_0,X_0+X_1)_{0,\infty,A;J}=(X_0,X_0+X_1)_{0,\infty,A;J;(1,\infty)}.
\end{equation}
}
\end{remark}

\smallskip
We shall also need the next two lemmas.

\begin{lemma}\label{01}
If $(X_0, X_1)$ is a compatible couple, $\theta \in [0,1]$, $1\le q\le \infty$, and 
$w \in  {\mcal W}(0,\infty),$ 
then 
\begin{equation}\label{X0X1}
(X_0, X_1)_{\theta, q, w; K}\equiv (X_1, X_0)_{1-\theta, q, w_1; K},
\end{equation}
where $w_1 \in  {\mcal W}(0,\infty)$ 
is given by
\begin{equation}\label{B1}
w_1(x):=w(1/x) \quad\text{for all\ }
x>0.
\end{equation}
\end{lemma}

\bpr
If $f \in X_0+X_1$, then 
\begin{equation}\label{L1}
\|f\|_{(X_0,X_1)_{\theta,q; w,K}}=\|t^{-\theta-1/q}w(t)K(f,t;X_0,X_1)\|_{q,(0,\infty)}.
\end{equation}
 Making use of \eqref{K} and a change of variables, we obtain that
\begin{align*}
{\rm RHS}\eqref {L1}&=\|t^{1-\theta-1/q}w(t)K(f,t^{-1};X_1,X_0)\|_{q,(0,\infty)}\\ \nonumber
&=\|s^{-(1-\theta)-1/q}w_1(s)K(f,s;X_1,X_0)\|_{q,(0,\infty)}\\ \nonumber
&=\|f\|_{(X_1,X_0)_{1-\theta,q; w_1,K}},
\end{align*}
and the result follows.
\epr

\smallskip
\begin{lemma}\label{10}
If $(X_0, X_1)$ is a compatible couple, $\theta \in [0,1]$, $1\le q\le \infty$, and 
$ w\in  {\mcal W}(0,\infty),$ 
then 
\begin{equation}\label{J0J1}
(X_0, X_1)_{\theta, q, w; J}\equiv(X_1, X_0)_{1-\theta, q, w_1; J},
\end{equation}
where $w_1 \in  {\mcal W}(0,\infty)$ 
is given by
\begin{equation}\label{A1}
w_1(x):=w(1/x) \quad \text{ for all} \ x>0.
\end{equation}
\end{lemma}

\bpr
If $f \in X_0+X_1$ and
\begin{equation}\label{1300}
f=\int_0^\infty u(t)\, \frac{dt}{t} \quad  {\rm (} {\text convergence\  in \ } X_0+X_1{\rm )},
\end{equation}
with $u(t) \in X_0\cap X_1$ for every $t>0,$ then also
$$
f=\int_0^\infty u_1(s)\, \frac{ds}{s} \quad \text{with } \ u_1(s)=u(1/s) \ \text{for every } s>0.
$$
Furthermore,   
\begin{equation}\label{L2}
\|f\|_{(X_0,X_1)_{\theta,q,w;J}}=\inf \|t^{-\theta-1/q}w(t)J(u(t),t;X_0,X_1)\|_{q,(0,\infty)},
\end{equation}
where the infimum extends over all representation \eqref{1300} of $f$. 
Making use of \eqref{J} and a~change of variables, we get 
\begin{align*}
\|t^{-\theta-1/q}w(t)J(u(t),t;X_0,X_1)\|_{q,(0,\infty)}
&=\|t^{1-\theta-1/q}w(t)J(u(t),t^{-1};X_1,X_0)\|_{q,(0,\infty)}\\ \nonumber
&=\|s^{-(1-\theta)-1/q}w_1(s)J(u_1(s),s;X_1,X_0)\|_{q,(0,\infty)},\\ \nonumber
\end{align*}
\vskip -0,5cm
\noindent
which, together with \eqref{L2}, implies that 
$$
\|f\|_{(X_0,X_1)_{\theta,q,w;J}}=\|f\|_{(X_1,X_0)_{1-\theta,q,w_1;J}}.
$$
Consequently, \eqref{J0J1} holds.
\epr

\smallskip
A description when the $K$-space $(X_0, X_1)_{\theta, q, B; K}$ with the limiting value 
$\theta=1$ and $B \in SV(0,\infty)$ coincides with the $J$-space $(X_0, X_1)_{\theta, q, A; J}$ with a  convenient $A \in SV(0,\infty)$ is given in the following assertion.

\begin{theorem} \label{ET1}
Let $(X_0, X_1)$ be a compatible couple and $1\le q <\infty$. If $B \in SV(0,\infty)$ satisfies
\begin{equation}\label{prop_slow_var_funct_B}
\int_0^x t^{-1}B^{\,q}(t)\, dt <\infty  \quad\text{for all \ }x>0, \qquad 
\int_0^\infty t^{-1}B^{\,q}(t)\, dt =\infty,
\end{equation}
and $A \in  SV(0,\infty)$ is such that 
\begin{equation}\label{def_slow_var_funct_A}
A(x) \approx B^{-q/q\,'}(x)\int_0^x t^{-1}B^{\,q}(t)\, dt  \quad\text{for a.a. \ }x>0,
\end{equation}
then 
\begin{equation}\label{1K_space=J_space1}
(X_0, X_1)_{1, q, B; K}=(X_0, X_1)_{1, q, A; J}.
\end{equation}
\end{theorem}

\bpr
By Lemma \ref{01},
\begin{equation}\label{aX0X1}
(X_0, X_1)_{1, q, B; K}\equiv (X_1, X_0)_{0, q, B_1; K},
\end{equation}
with $B_1$ given by $B_1(x):=B(1/x)$ for all $x>0$. 

Using assumptions \eqref{prop_slow_var_funct_B} and a change of variables, we can see that 
$$
\int_x^\infty t^{-1}B_1^{\,q}(t)\, dt <\infty  \quad\text{for all \ }x>0, \qquad 
\int_0^\infty t^{-1}B_1^{\,q}(t)\, dt =\infty.
$$%
Thus, we can apply Theorem \ref{equivalence theorem1} to arrive at
\begin{equation}\label{11K_space=J_space1}
(X_1, X_0)_{0, q, B_1; K}=(X_1, X_0)_{0, q, \tilde{A}; J}\,,
\end{equation}
with $\tilde{A} \in  SV(0,\infty)$ satisfying 
\begin{align}\label{1def_slow_var_funct_a1}
\tilde{A}(x)&\approx B_1^{-q/q\,'}(x)\int_x^\infty t^{-1}B_1^{\,q}(t)\, dt\\ \nonumber
&=B^{-q/q\,'}(1/x)\int_0^{1/x} t^{-1}B^{\,q}(t)\, dt 
 \quad\text{for a.a.\ }x>0\,. \nonumber
\end{align}
Consequently, by Lemma \ref{10},
\begin{equation}\label{111K_space=J_space1}
(X_1, X_0)_{0, q, \tilde{A}; J}=(X_0, X_1)_{1, q, A; J},
\end{equation}
with $A \in  SV(0,\infty)$ given by $A(x)=\tilde{A}(1/x)$ for all $x>0.$ 
Together with \eqref{1def_slow_var_funct_a1}, this implies that the function $A$ satisfies~\eqref{def_slow_var_funct_A}.

On using \eqref{aX0X1}, \eqref{11K_space=J_space1}, and \eqref{111K_space=J_space1}, we also obtain 
\eqref{1K_space=J_space1}.
\epr

\smallskip
There is the following counterpart of the previous result.

\begin{theorem} \label{ET1.1}
Let $ \left(X_0,X_1\right)$  be a compatible couple and $1<q \leq \infty$. If $A\in SV(0,\infty)$ satisfies
\begin{equation}\label{c1.1}
{\int_{x}^{\infty} t^{-1}A^{-q\,'}(t)\,dt}< \infty \hspace{0.2cm}\text{  for all $x>0$, }\hspace{0.5cm}{\int_0^{\infty} t^{-1}A^{-q\,'}(t)\,dt}=\infty,
\end{equation}   and $B\in SV(0,\infty)$ is such that
\begin{equation}\label{a.1}
B(x)\approx A^{-{q\,'\!}/{q}}(x)\left(\int\limits_{x}^{\infty}t^{-1}A^{-q\,'}(t)\,dt\right)^{-1}\quad \text{for a.a. }x>0,
\end{equation}
then 
\begin{equation}\label{01K_space=J_space1}
(X_0,X_1)_{1,q,A;J}=(X_0,X_1)_{1,q,B;K}\,.
\end{equation}
\end{theorem}

\bpr
By Lemma \ref{10},
\begin{equation}\label{aaX0X1}
(X_0, X_1)_{1, q, A; J}\equiv (X_1, X_0)_{0, q, A_1; J},
\end{equation}
with $A_1$ given by $A_1(x):=A(1/x)$ for all $x>0$. 

Using assumptions \eqref{c1.1} and a change of variables, we can see that 
$$
\int_0^x t^{-1}A_1^{q\,'}(t)\, dt <\infty  \quad\text{for all \ }x>0, \qquad 
\int_0^\infty t^{-1}A_1^{q\,'}(t)\, dt =\infty.
$$%
Thus, applying Theorem \ref{equivalence theorem}, we arrive at
\begin{equation}\label{11J_space=K_space1}
(X_1, X_0)_{0, q, A_1; J}=(X_1, X_0)_{0, q, \tilde{B}; K}\,,
\end{equation}
with $\tilde{B} \in  SV(0,\infty)$ satisfying
\begin{align}\label{11def_slow_var_funct_a1}
\tilde{B}(x)&:=A_1^{-q\,'\!/q}(x)\Big(\int_0^x t^{-1}A_1^{-q\,'}(t)\, dt\Big)^{-1}\\ \nonumber
&=A^{-q\,'\!/q}(1/x)\Big(\int_{1/x}^\infty t^{-1}A^{-q\,'}(t)\, dt \Big)^{-1}
 \quad\text{for a.a. \ }x>0\,. \nonumber
\end{align}
Consequently, by Lemma \ref{01},
\begin{equation}\label{111K_space=J_space11}
(X_1, X_0)_{0, q, \tilde{B}; K}=(X_0, X_1)_{1, q, B; K},
\end{equation}
with $B \in  SV(0,\infty)$ given by $B(x)=\tilde{B}(1/x)$ for all $x>0.$ 
Together with \eqref{11def_slow_var_funct_a1}, this implies that the function $B$ satisfies~\eqref{a.1}.

On using \eqref{aaX0X1}, \eqref{11J_space=K_space1}, and \eqref{111K_space=J_space11}, we also obtain 
\eqref{01K_space=J_space1}.
\epr

\smallskip
We continue with 
 discrete characterizations of spaces $(X_0,X_1)_{\theta,q,b; K}$ and $(X_0,X_1)_{\theta,q,b; J}$. To this end we shall use the following assertion.

\begin{lemma}\label{mind}
If $\theta \in [0,1],$ and $b\in SV(0,\infty),$ then there are positive constants $c_1, c_2$ such that
$$
\label{mindc}
\frac{c_1}{2^{\theta}} \,2^{-m\theta}b(2^m)\leq t^{-\theta}b(t) \leq c_2\, 2^{-m\theta}b(2^m)\quad \text{for all } \ \theta \in [0,1],
\ m\in \mathbb{Z} \ \text{and } \ t\in[2^m, 2^{m+1}).
$$
\end{lemma}

\bpr 
Since the function $t\mapsto tb(t),\  t\in (0,\infty),$ is equivalent to a non-decreasing function, there is a positive constant $k_1$ such that
$$
2^mb(2^m)\le k_1 tb(t)\quad \text{for all } m\in \mathbb{Z} \ \text{and } t\in[2^m, 2^{m+1}),
$$
which implies that
\begin{equation}\label{1mindI1}
\frac{b(2^m)}{k_1}\le \frac{tb(t)}{2^m}\le \frac{2^{m+1}b(t)}{2^m}=2b(t)\quad \text{for all } m\in \mathbb{Z} \ \text{and } t\in[2^m, 2^{m+1}).
\end{equation}
Analogously, since the function $t\mapsto \frac{b(t)}{t}, \ t\in (0,\infty),$ is equivalent to a non-increasing function, there is a~positive constant $k_2$ such that
$$
b(t)=t\frac{b(t)}{t}\le t\,k_2\frac{b(2^m)}{2^m}\le k_2\,\frac{2^{m+1}}{2^m} b(2^m)=2k_2\,b(2^m)\quad \text{for all } m\in \mathbb{Z} \ \text{and } t\in[2^m, 2^{m+1}),
$$
which, together with \eqref{1mindI1}, implies that 
\begin{equation}\label{2mindI1}
\frac{b(2^m)}{2k_1}\le b(t)\le 2k_2 \,b(2^m)\quad \text{for all } m\in \mathbb{Z} \ \text{and } t\in[2^m, 2^{m+1}).
\end{equation}
Moreover, the function $t\mapsto t^{-\theta},\  t\in (0,\infty),$ with $\theta \in [0, 1],$ is non-increasing. Consequently,
\begin{equation}\label{3mindI1}
\big(2^{m+1}\big)^{-\theta} \le t^{-\theta} \le \big(2^{m}\big)^{-\theta} \quad \text{for all } \ m\in \mathbb{Z} \ \text{and } \ t\in[2^m, 2^{m+1}).
\end{equation}
Using estimates \eqref{2mindI1} and \eqref{3mindI1}, we arrive at
$$
\big(2^{m+1}\big)^{-\theta}\, \frac{b(2^m)}{2k_1}\le t^{-\theta}b(t)\le \big(2^{m}\big)^{-\theta}\,2k_2\, b(2^m) \quad \text{for all } 
\ \theta \in [0,1],\ m\in \mathbb{Z} \ \text{and } \ t\in[2^m, 2^{m+1})\,,
$$
and the result 
follows (with $c_1=1/(2k_1)$ and $c_2=2k_2$).
\epr

\smallskip
\begin{lemma}\label{KD}
If $(X_0,X_1)$ is a compatible couple, $\theta\in [0,1], 1\leq q\leq \infty,$ and $ b\in SV(0,\infty),$ then 
there are positive constants $c_1, c_2$ such that
\begin{equation}\label{KDE}
\frac{c_1}{2^{\theta}}(\ln 2)^{1/q}\,\lVert \{K(f,2^m)\}_{m \in \mathbb{Z}} \rVert_{\lambda_{\theta,q,b}}\le \lVert f\rVert_{\theta,q,b;K} 
\le 2c_2(\ln 2)^{1/q}\,\lVert \{K(f,2^m)\}_{m \in \mathbb{Z}} \rVert_{\lambda_{\theta,q,b}} 
\end{equation}
for all $f \in X_0 +X_1$, 
$\theta\in [0,1]$, and $1\leq q\leq \infty.$

In particular, given $f \in X_0 +X_1$, then $f\in (X_0,X_1)_{\theta,q,b;K}$ if and only if $\{K(f,2^m)\}_{m\in \mathbb{Z}} \in \lambda_{\theta,q,b}$ and
\begin{equation}\label{KDE1}
\lVert f\rVert_{\theta,q,b;K} \approx 
\lVert \{K(f,2^m)\}_{m \in \mathbb{Z}} \rVert_{\lambda_{\theta,q,b}}
\end{equation}
for all $f \in (X_0,X_1)_{\theta,q,b; K}$,  
$\theta\in [0,1]$, and $1\leq q\leq \infty.$
\end{lemma}
\bpr 
Let $f \in X_0+X_1.$ 
Using properties of the $K$-functional (cf. Lemma \ref{about K}), we obtain
$$
K(f,2^m) \le K(f,t) \le 2K(f,2^m)\quad\text{for all } m\in \mathbb{Z} \ \text{and } t\in[2^m, 2^{m+1}).
$$
Together with Lemma \ref{mind}, this implies that 
\begin{equation}\label{2121}
\frac{c_1}{2^{\theta}}\,2^{-m\theta}b(2^m)K(f,2^m)\le t^{-\theta}b(t)K(f,t) \le 2c_2\,2^{-m\theta}b(2^m)K(f,2^m) 
\end{equation}
for all $\theta \in [0,1],\ m\in \mathbb{Z} \ \text{and } \ t\in[2^m, 2^{m+1})\,.$

Hence, if $1\le q<\infty$, then 
\begin{align}\label{2122}
\Big(\frac{c_1}{2^{\theta}}\,2^{-m\theta}b(2^m)K(f,2^m)\Big)^q\,\ln 2 
&\le \int_{2^m}^{2^{m+1}}(t^{-\theta}b(t)K(f,t))^q\frac{dt}{t}\\ \nonumber
&\le \Big(2c_2\,2^{-m\theta}b(2^m)K(f,2^m)\Big)^q \,\ln 2 ,
\end{align}
which implies that 
\begin{align*}
\frac{c_1}{2^{\theta}}\,(\ln 2)^{1/q}\,\Big( \sum_{m\in\mathbb{Z}}\big(2^{-m\theta}b(2^m)&K(f,2^{m})\big)^q\Big)^{1/q}\\ \nonumber
&\le \lVert f\rVert_{\theta,q,b;K} \\ \nonumber
&\le 2c_2\,(\ln 2)^{1/q}\,\Big( \sum_{m\in\mathbb{Z}}\big(2^{-m\theta}b(2^m)K(f,2^{m})\big)^q\Big)^{1/q},
\end{align*}
i.e,
$$
\frac{c_1}{2^{\theta}}\,(\ln 2)^{1/q}\,\lVert \{K(f,2^m)\}_{m \in \mathbb{Z}} \rVert_{\lambda_{\theta,q,b}} 
\le \lVert f\rVert_{\theta,q,b;K} 
\le 2c_2\,(\ln 2)^{1/q}\,\lVert \{K(f,2^m)\}_{m \in \mathbb{Z}} \rVert_{\lambda_{\theta,q,b}} 
$$
for all $f \in X_0+X_1$, $\theta \in [0,1]$ and $1\le q< \infty,$ which is \eqref{KDE} with $1\le q<\infty.$

Furthermore, if $q=\infty,$ then using \eqref{2121}, we obtain  
$$
\frac{c_1}{2^{\theta}}\,2^{-m\theta}b(2^m)K(f,2^m)\le \esssup_{t\in [2^m, 2^{m+1})} t^{-\theta}b(t)K(f,t)
\le 2c_2\,\,2^{-m\theta}b(2^m)K(f,2^m)
\quad \text{for all } m\in \mathbb{Z}, 
$$
which immediately gives 
$$
\frac{c_1}{2^{\theta}}\lVert \{K(f,2^m)\}_{m \in \mathbb{Z}} \rVert_{\lambda_{\theta,q,b}} 
\le \lVert f\rVert_{\theta,q,b;K} 
\le 2c_2 \lVert \{K(f,2^m)\}_{m \in \mathbb{Z}} \rVert_{\lambda_{\theta,q,b}} 
$$
for all $f \in X_0+X_1, \  \theta \in [0,1]$, and $q=\infty$, 
which is \eqref{KDE} with $q=\infty.$
\epr

\smallskip
 \begin{lemma}\label{JD} Let $(X_0,X_1)$ be a compatible couple, $\theta \in [0,1], 1\leq q\leq \infty,$ and $b\in SV(0,\infty)$. Then
$f\in (X_0,X_1)_{\theta,q,b;J}$ if and only if there exist $u_m\in X_0\cap X_1$ for all $m\in\mathbb{Z}$ such that
\begin{equation}\label{JDBI}
f=\sum_{m\in\mathbb{Z}} u_m\hspace{1cm}\text{{\rm (}convergence in $X_0+X_1${\rm )}}
\end{equation}
and
\begin{equation}\label{JDC}
\{J(u_m,2^m)\}_{m\in \mathbb{Z}}\in \lambda_{\theta,q,b}.
\end{equation}
Moreover, there are positive constants $c_1, c_2$ such that
\begin{equation}\label{JDS}
\frac{1}{2c_2}(\ln 2)^{1/q\,'}\,\lVert f\rVert_{\theta,q,b;J}\le\inf
\,\lVert \{J(u_m,2^m)\}_{m\in \mathbb{Z}}\rVert_{\lambda_{\theta,q,b}}\le\frac{2^\theta}{c_1}(\ln 2)^{1/q\,'} \lVert f\rVert_{\theta,q,b;J}
\end{equation} 
for all $f \in (X_0,X_1)_{\theta,q,b; J},\, \theta\in [0,1]$, and $1\leq q\leq \infty,$ 
where the infimum is extended over all sequences $\{u_m\}_{m\in \mathbb{Z}}$ satisfying \eqref{JDBI} and \eqref{JDC}.
\end{lemma}
\bpr
Let  $f\in (X_0,X_1)_{\theta,q,b;J}$ and 
\begin{equation}\label{JDBI2}
f=\int_0^{\infty}u(t)\frac{dt}{t}\qquad\text{(convergence in $X_0+ X_1$)},
\end{equation}
where $u :(0,\infty)\rightarrow X_0 \cap X_1$ is a strongly measurable function.  
Putting $u_m=\int_{2^m}^{2^{m+1}}u(t)\,dt/t$ for all $m\in\mathbb{Z},$ we see that \eqref{JDBI} holds. Since  
\begin{equation}\label{subadJ}
J(u_m, 2^m)= J\Big(\int_{2^m}^{2^{m+1}} u(t)\frac{dt}{t}, 2^m\Big)\le \int_{2^m}^{2^{m+1}} J(u(t), 2^m)\frac{dt}{t}
\le \int_{2^m}^{2^{m+1}} J(u(t),t)\frac{dt}{t}
\end{equation}
for all $m\in\mathbb{Z}$ and $t\in <2^m, 2^{m+1})$, we get
\begin{equation}\label{JDE1}
\lVert \{J(u_m,2^m)\}_{m\in\mathbb{Z}}\rVert_{\lambda_{\theta,q,b}} 
\le\left\lVert \left\lbrace\int_{2^m}^{2^{m+1}}J(u(t),t)\frac{dt}{t}\right\rbrace_{m\in\mathbb{Z}}\right\rVert_{\lambda_{\theta,q,b}}.
\end{equation}

If $1\le q<\infty,$ then, by Lemma \ref{mind}, 
\begin{align}\label{JDE1.1}
{\rm RHS\eqref{JDE1}}&=\left(\sum_{m\in\mathbb{Z}}\left(2^{-m\theta}b(2^m)\int_{2^m}^{2^{m+1}}J(u(t),t)\frac{dt}{t}\right)^q\right)^{1/q}\\
&\le \frac{2^\theta}{c_1} \left(\sum_{m\in\mathbb{Z}}\left(\int_{2^m}^{2^{m+1}}t^{-\theta}b(t)J(u(t),t)\frac{dt}{t}\right)^q\right)^{1/q}.\nonumber
\end{align}

Moreover, by the Jensen inequality, 
$$
\left(\int_{2^m}^{2^{m+1}}t^{-\theta}b(t)J(u(t),t)\frac{dt}{t \ln 2}\right)^q 
\le \int_{2^m}^{2^{m+1}}(t^{-\theta}b(t)J(u(t),t))^q\frac{dt}{t \ln 2} \quad \text{for all}\  m\in\mathbb{Z}.
$$
Together with \eqref{JDE1.1} and \eqref{JDE1}, this gives 
$$
\lVert \{J(u_m,2^m)\}_{m\in\mathbb{Z}}\rVert_{\lambda_{\theta,q,b}}
\le\frac{2^\theta}{c_1}(\ln 2)^{1/q\,'} \|t^{-\theta-1/q}b(t)J(u(t),t)\|_{q, (0,\infty)},
$$
which implies \eqref{JDC}. Thus, taking the infimum over all sequences $\{u_m\}_{m\in \mathbb{Z}}$ satisfying \eqref{JDBI} and \eqref{JDC}, 
we get the second inequality in \eqref{JDS}  
for all $f \in (X_0,X_1)_{\theta,q,b;J},\, \theta\in [0,1]$, and $1\leq q < \infty.$


If $q=\infty$, then, using \eqref{JDE1} and Lemma \ref{mind}, we obtain 
$$
{\rm RHS\eqref{JDE1}}=\sup_{m\in\mathbb{Z}}2^{-m\theta}b(2^m)\int_{2^m}^{2^{m+1}}J(u(t),t)\frac{dt}{t}\\
\le \frac{2^\theta}{c_1} \sup_{m\in\mathbb{Z}}\,\int_{2^m}^{2^{m+1}}t^{-\theta}b(t)J(u(t),t)\frac{dt}{t}.\nonumber
$$
Since (cf. \cite[Lemma 5.5, p. 47]{OK90:HTI})
$$
t^{-\theta}b(t)J(u(t),t) \le \esssup_{s\in [2^m, 2^{m+1})} s^{-\theta}b(s)J(u(s),s)\quad \text{for all}\ m\in\mathbb{Z}\ \text{and}\ \text{for a.a.\ }t\in [2^m,2^{m+1}),
$$
we arrive at
\begin{align}\label{nic}
{\rm RHS\eqref{JDE1}}&=\sup_{m\in\mathbb{Z}}2^{-m\theta}b(2^m)\int_{2^m}^{2^{m+1}}J(u(t),t)\frac{dt}{t}\\
&\le \frac{2^\theta}{c_1} (\ln 2) \,\sup_{m\in\mathbb{Z}}\,\esssup_{s\in <2^m, 2^{m+1})} s^{-\theta}b(s)J(u(s),s)\nonumber\\
&=\frac{2^\theta}{c_1} (\ln 2)\,\|s^{-\theta-1/q}b(s)J(u(s),s)\|_{q, (0,\infty)}.\nonumber
\end{align}
Together with \eqref{JDE1}, this gives 
$$
\lVert \{J(u_m,2^m)\}_{m\in\mathbb{Z}}\rVert_{\lambda_{\theta,q,b}}
\le\frac{2^\theta}{c_1}(\ln 2)\, \|t^{-\theta-1/q}b(s)J(u(s),s)\|_{q, (0,\infty)},
$$
which implies that  
the second inequality in \eqref{JDS} holds for all $f \in X_0,X_1)_{\theta,q,b; J},\, \theta\in [0,1]$, and $q=\infty.$

Conversely, assume that
\begin{equation}\label{JDBI3}
f=\sum_{m\in\mathbb{Z}} u_m\quad\text{(convergence in $X_0+X_1$)},
\end{equation}
where $u_m\in X_0\cap X_1$ for all $m\in\mathbb{Z}$, and 
\begin{equation}\label{JDCr}
\{J(u_m,2^m)\}_{m\in\mathbb{Z}}\in {\lambda_{\theta,q,b}}.
\end{equation}
Putting 
\begin{equation}\label{ufce}
u(t)=\frac{u_m}{\ln 2} \quad \text{if} \ \ 2^m\leq t<2^{m+1}\  \text{and} \ \  m\in\mathbb{Z},
\end{equation}
we get $u:(0,\infty)\rightarrow X_0\cap X_1$ and 
$$f=\sum_{m\in\mathbb{Z}}u_m=\sum_{m\in\mathbb{Z}}\int_{2^m}^{2^{m+1}}\frac{u_m}{\log{2}}\frac{dt}{t}=\int_0^{\infty}u(t)\frac{dt}{t}.$$
Thus,
\begin{equation}\label{JDE2.1}
\lVert f\rVert_{\theta,q,b;J}\leq \lVert t^{-\theta-{1}/{q}}b(t)J(u(t),t)\rVert_{q,(0,\infty)}.
\end{equation}

If $1\leq q<\infty$, then using \eqref{ufce}, 
we obtain 
\begin{align*}
{\rm RHS\eqref{JDE2.1}}&=\left(\sum_{m\in\mathbb{Z}}\int_{2^m}^{2^{m+1}}\left(t^{-\theta}b(t)J\left(\frac{u_m}{\ln 2},t\right)\right)^q\frac{dt}{t}\right)^{1/q}\\
&=\frac{1}{\ln 2}\left(\sum_{m\in\mathbb{Z}}\int_{2^m}^{2^{m+1}}\left(t^{-\theta}b(t)J\left(u_m,t\right)\right)^q\frac{dt}{t}\right)^{1/q}.
\end{align*}
Together with the estimate 
\begin{equation}\label{Jfun}
J(u_m, t)\le 2\,J(u_m, 2^m)\quad \text{for all} \ m \in \mathbb{Z}\ \text{and} \  t \in [2^m, 2^{m+1}),
\end{equation}
and Lemma \ref{mind}, this implies that
\begin{align}\label{pravastr}
{\rm RHS\eqref{JDE2.1}}&\le\frac{2}{\ln 2}\,c_2\left(\sum_{m\in\mathbb{Z}}\int_{2^m}^{2^{m+1}}\left(2^{-m\theta}b(2^m)J\left(u_m,2^m\right)\right)^q\frac{dt}{t}\right)^{1/q}\\
&=\frac{2}{\ln 2}\,c_2\,(\ln 2)^{1/q}\left(\sum_{m\in\mathbb{Z}}\left(2^{-m\theta}b(2^m)J\left(u_m,2^m\right)\right)^q\right)^{1/q}\nonumber\\
&=2c_2 (\ln 2)^{-1/q\,'}\lVert \{J(u_m,2^m)\}_{m\in\mathbb{Z}}\rVert_{\lambda_{\theta,q,b}}.\nonumber
\end{align}
Combining estimates \eqref{JDE2.1} and \eqref{pravastr}, we arrive at 
 $$
\frac{1}{2c_2}(\ln 2)^{1/q\,'}\,\lVert f\rVert_{\theta,q,b;J}\
\le \lVert \{J(u_m,2^m)\}_{m\in\mathbb{Z}}\rVert_{\lambda_{\theta,q,b}},
$$
and taking the infimum over all sequences $\{u_m\}_{m\in \mathbb{Z}}$ satisfying \eqref{JDBI3} and \eqref{JDCr},
we get the first inequality in \eqref{JDS}.     
     
If $q=\infty$, then, by \eqref{ufce},
\begin{align*}
{\rm RHS\eqref{JDE2.1}}&=\sup_{m\in\mathbb{Z}}\,\esssup_{s\in [	2^m, 2^{m+1})} s^{-\theta}b(s)J(u(s),s)\\
&=\frac{1}{\ln 2}\,\sup_{m\in\mathbb{Z}}\,\esssup_{s\in [2^m, 2^{m+1})} s^{-\theta}b(s)J(u_m,s).
\end{align*}
Using estimate \eqref{Jfun} and Lemma \ref{mind}, we obtain
\begin{align}\label{lab}
{\rm RHS\eqref{JDE2.1}}&\le \frac{2}{\ln 2}c_2\,\sup_{m\in\mathbb{Z}}\,\esssup_{s\in [2^m, 2^{m+1})} 2^{-m\theta}b(2^m)J(u_m, 2^m)\\
&=\frac{2}{\ln 2}c_2\,\sup_{m\in\mathbb{Z}} 2^{-m\theta}b(2^m)J(u_m, 2^m)\nonumber\\
&=\frac{2}{\ln 2}c_2\,\lVert \{J(u_m,2^m)\}_{m\in\mathbb{Z}}\rVert_{\lambda_{\theta,q,b}}.\nonumber
\end{align}
Combining estimates \eqref{JDE2.1} and \eqref{lab}, we arrive at 
 $$
\frac{1}{2c_2}(\ln 2)\,\lVert f\rVert_{\theta,q,b;J}\
\le \lVert \{J(u_m,2^m)\}_{m\in\mathbb{Z}}\rVert_{\lambda_{\theta,q,b}},
$$
and taking the infimum over all sequences $\{u_m\}_{m\in \mathbb{Z}}$ satisfying \eqref{JDBI3} and \eqref{JDCr},
we get the first inequality in \eqref{JDS}. 
\epr

\begin{remark}\label{cc} {\rm Note that the constants $c_1$ and $c_2$ in Lemmas \ref{KD} and \ref{JD}    are the same as those in Lemma \ref{mind}   and thus they depend only on the function $b\in SV(0,\infty)$. 
}
\end{remark}

\smallskip
We shall also need the following assertion.

\begin{lemma}\label{JDchara} 
Let $\theta \in [0, 1]$, $1\le q \le \infty$, and $a \in  SV(0,\infty)$. Then  
\begin{equation}\label{DC1}
\left\|t^{\theta-1/q\,'}\,a^{-1}(t) \min\left\{1, t^{-1}\right\}\right\|_{q\,'\!, (0,\infty)}<\infty,
\end{equation}
if and only if
\begin{equation}\label{307bd}
\{\min \{1, 2^m\}\}_{m\in \mathbb{Z}}\in \lambda_{1-\theta,q\,'\!,a^{-1}}. 
\end{equation}
\end{lemma}

\bpr
Note that 
\begin{equation}\label{odhad}
\frac{1}{2}\,\min\{1,t\} \le \frac{1}{2}\,\min \{1, 2^{m+1}\}\le \min \{1, 2^m\}\le \min\{1,t\}
\end{equation}
for all $m\in \mathbb{Z}$ and $t\in [2^m, 2^{m+1})$.

If $1<q\le\infty$, then, using Lemma \ref{mind} and  \eqref{odhad}, we obtain 
\begin{align*}
\Big\|\{\min \{1, 2^m\}\}_{m\in \mathbb{Z}}\Big\|_{\lambda_{1-\theta,q\,'\!,a^{-1}}}&=
\Big(\sum_{m\in \mathbb{Z}} \big(2^{-m(1-\theta)} a^{-1}(2^m) \min \{1, 2^m\}\big)^{q\,'}\Big)^{1/q\,'}\\
&\approx\Big(\sum_{m\in \mathbb{Z}} \big(2^{-m(1-\theta)} a^{-1}(2^m) \min \{1, 2^m\}\big)^{q\,'}\int_{2^m}^{2^{m+1}}\frac{dt}{t}\Big)^{1/q\,'}\\
&\approx\Big(\sum_{m\in \mathbb{Z}} \int_{2^m}^{2^{m+1}}\big(t^{-(1-\theta)} a^{-1}(t) \min \{1, t\}\big)^{q\,'}\frac{dt}{t}\Big)^{1/q\,'}\\
&=\Big(\int_0^\infty\big(t^{\theta} a^{-1}(t) \min \{1, t^{-1}\}\big)^{q\,'}\frac{dt}{t}\Big)^{1/q\,'}\\
&=\big\|t^{\theta-1/q\,'} a^{-1}(t) \min \{1, t^{-1}\}\big\|_{q\,'\!, \,(0,\infty)},
\end{align*}
and the result follows.

Assume now that $q=1$. Then, by Lemma \ref{mind} and  \eqref{odhad},
\begin{align*}
\Big\|\{\min \{1, 2^m\}\}_{m\in \mathbb{Z}}\Big\|_{\lambda_{1-\theta,q'\!,a^{-1}}}&=
\sup_{m\in \mathbb{Z}}\, 2^{-m(1-\theta)} a^{-1}(2^m) \min \{1, 2^m\}\\
&\approx \sup_{m\in \mathbb{Z}}\ \esssup_{s\in [2^m, 2^{m+1})}t^{-(1-\theta)} a^{-1}(t) \min \{1, t\}\\
&=\big\|t^{-(1-\theta)} a^{-1}(t) \min \{1, t\}\big\|_{q\,'\!, \,(0,\infty)}\\
&=\big\|t^{\theta-1/q\,'} a^{-1}(t) \min \{1, t^{-1}\}\big\|_{q\,'\!, \,(0,\infty)},
\end{align*}
and the result follows.
\epr

\smallskip
\begin{theorem}\label{E1}
Let $(X_0,X_1)$ be a compatible couple, $\theta \in [0, 1]$, $1\leq q<\infty$, and let $X_0\cap X_1$ be dense in $X_0$ and $X_1$. If $a\in SV(0,\infty)$ satisfies condition \eqref{DC1}, 
then 
\begin{equation}\label{embdg1}
(X_0,X_1)'_{\theta,q,a;J}\hookrightarrow (X_1',X_0')_{1-\theta,q\,'\!,a^{-1};K}.
\end{equation}
\end{theorem}
\begin{proof}
By Theorem \ref{305},  condition \eqref{DC1} guarantees that the space $(X_0,X_1)_{\theta,q,a;J}$ is a~Banach space, which is an intermediate space between $X_0$ and $X_1$. In particular, 
\begin{equation}\label{E1.1.1}
X_0\cap X_1\hookrightarrow (X_0,X_1)_{\theta,q,a;J}.
\end{equation}
 Moreover, by \eqref{99}, 
\begin{equation}\label{E1.1.2}
(X_0\cap X_1)'=X_0'+X_1'.
\end{equation}
Let $0\neq f'\in (X_0,X_1)'_{\theta,q,a;J}$. By \eqref{E1.1.1} and \eqref{E1.1.2},
\begin{equation}\label{E1.1.3}
f'\in X_0'+X_1'.
\end{equation}
Since $K(f',2^m;X_0',X_1')$ is the equivalent norm on $X_0'+X_1'$ for every $m\in\mathbb{Z}$, from \eqref{E1.1.3} we see that 
\begin{equation}\label{E1.1,4}
0<K(f', 2^m; X_0',X_1')<\infty.
\end{equation}
Given $\varepsilon>0$ and $m\in\mathbb{Z}$. Then \eqref{010} implies that there is $w_m\in X_0\cap X_1$ such that
\begin{equation}\label{epsilon1}
K(f',2^{-m};X_0',X_1')-\varepsilon\min{\{1,2^{-m}\}}\leq \frac{\langle f',w_m\rangle}{J(w_m,2^m)}\,.
\end{equation}
(Note that, by \eqref{010} and \eqref{E1.1,4}, $J(w_m,2^m)<\infty$). Choose a sequence of numbers $\tau:=\{\tau_m\}_{m \in \mathbb{Z}}\in \lambda_{\theta,q,a},$ and put 
\begin{equation}\label{E1BS}
f_{\tau}:=\sum_{m\in\mathbb{Z}}u_m\qquad\text{(convergence in $X_0+X_1$)}\,,
\end{equation}
where
\begin{equation}\label{BSUM}
u_m:=\frac{\tau_m w_m}{J(w_m,2^m)}\in X_0\cap X_1,\qquad\text{for all }m\in\mathbb{Z}\,.
\end{equation}
Then
\begin{align*}
\lVert \{J(u_m,2^m)\}_{m\in\mathbb{Z}}\rVert_{\lambda_{\theta,q,a}}&=\left(\sum_{m\in\mathbb{Z}}\left(2^{-m\theta}a(2^m)J\left(\frac{\tau_m w_m}{J(w_m,2^m)},2^m\right)\right)^q \right)^{{1}/{q}}\\
&=\left(\sum_{m\in\mathbb{Z}}\left(2^{-m\theta}a(2^m)\frac{\tau_m}{J(w_m,2^m)}J\left(w_m,2^m\right)\right)^q \right)^{{1}/{q}}\\
&=\lVert \{\tau_m\}_{m\in\mathbb{Z}}\rVert_{\lambda_{\theta,q,a}}<\infty.
\end{align*}
Thus, by Lemma \ref{JD}, 
\begin{equation}\label{666}
f_{\tau}\in (X_0,X_1)_{\theta,q,a;J}\qquad \text{and }\qquad  
\lVert f_{\tau}\rVert_{\theta,q,a;J}\lesssim \lVert \{\tau_{m}\}_{m \in \mathbb{Z}}\rVert_{\lambda_{\theta,q,a}}\,.
\end{equation}
Consequently,
\begin{align}\label{thus1.1}
\mid \langle f',f_{\tau}\rangle \mid &\leq \lVert f'\rVert_{(X_0,X_1)'_{\theta,q,a;J}}.\,\lVert f_{\tau}\rVert_{\theta,q,a;J}\\
&\lesssim \lVert f'\rVert_{(X_0,X_1)'_{\theta,q,a;J}}.\,\lVert \{\tau_m\}_{m \in \mathbb{Z}}\rVert_{\lambda_{\theta,q,a}}.\nonumber
\end{align}
Moreover, by \eqref{E1BS}, \eqref{BSUM}, and \eqref{epsilon1}, 
\begin{align}\label{grtrE}
\langle f',f_{\tau}\rangle &= \sum_{m\in\mathbb{Z}} \tau_m\,\frac{\langle f',w_m\rangle}{J(w_m,2^m)}\\
&\geq \sum_{m\in\mathbb{Z}}\tau_m\,(K(f',2^{-m};X_0',X_1')-\varepsilon\min\{1,2^{-m}\}).\nonumber
\end{align}
Together with \eqref{thus1.1}, \eqref{666}, and the fact that
$$
K(f',2^{-m};X_0',X_1')=2^{-m}K(f',2^m;X_1',X_0'),
$$
this implies that
$$
\sum_{m\in\mathbb{Z}} 2^{-m}\tau_m(K(f',2^m;X_1',X_0')-\varepsilon\min{\{1,2^m\}})\lesssim \lVert f'\rVert_{(X_0,X_1)_{\theta,q,a;J}'}\lVert \{\tau_m\}_{m \in \mathbb{Z}}\rVert_{\lambda_{\theta,q,a}},
$$
i.e.,
\begin{align}\label{epsilon2}
\sum_{m\in\mathbb{Z}} 2^{-m}\tau_m K(f',2^m;X_1',X_0')&\lesssim\lVert f'\rVert_{(X_0,X_1)_{\theta,q,a;J}'}\lVert \{\tau_m\}_{m \in \mathbb{Z}}\rVert_{\lambda_{\theta,q,a}}\\
&\qquad+ \varepsilon\sum_{m\in\mathbb{Z}} 2^{-m}\tau_m \min{\{1,2^m\}}. \nonumber
\end{align}
By Remark \ref{dual of lamba space},
$$
\sum_{m\in\mathbb{Z}} 2^{-m}{\tau_m}\min{\{1,2^m\}}\leq \lVert \{\tau_m\}_{m \in \mathbb{Z}}\rVert_{\lambda_{\theta,q,a}}\,
\lVert \{\min{\{1,2^m\}}\}_{m \in \mathbb{Z}}\rVert_{\lambda_{1-\theta,q'\!,a^{-1}}}.
$$
Since, by our assumption \eqref{DC1} and Lemma \ref{JDchara},
$$
c:=\lVert \{\min{\{1,2^m\}}\}_{m \in \mathbb{Z}}\rVert_{\lambda_{1-\theta,q',a^{-1}}}<\infty.
$$
we get from \eqref{epsilon2} that
$$
\sum_{m\in\mathbb{Z}} 2^{-m}\tau_m K(f',2^m;X_1',X_0')\lesssim (\lVert f'\rVert_{(X_0,X_1)_{\theta,q,a;J}'}+c \,\varepsilon)\,\lVert \{\tau_m\}_{m \in \mathbb{Z}}\rVert_{\lambda_{\theta,q,a}}\,,
$$
for all positive $\varepsilon$, all sequences $\tau:=\{\tau_m\}_{m \in \mathbb{Z}}\in \lambda_{\theta,q,a}$ and all $f'\in (X_0,X_1)'_{\theta,q,a;J}$.
Together with Remark \ref {dual of lamba space}, this implies that
\begin{equation}\label{E1KD}
\lVert \{K(f',2^m;X_1',X_0')\}\rVert_{\lambda_{1-\theta,q',a^{-1}}}\lesssim \lVert f'\rVert_{(X_0,X_1)'_{\theta,q,a;J}}
\quad \text{all }\  f'\in (X_0,X_1)'_{\theta,q,a;J}.
\end{equation}
By Lemma \ref{KD}, 
$$
{\rm LHS}\eqref{E1KD} \approx \lVert f'\rVert_{(X_1',X_0')_{1-\theta,q'\!, a^{-1};K}}\,.
$$
This estimate and  \eqref{E1KD} show that embedding \eqref{embdg1} holds.
\epr

\smallskip
\begin{theorem}\label{E2}
Let $(X_0,X_1)$ be a compatible couple, $\theta\in [0, 1]$, $1\leq q<\infty$, and let $X_0\cap X_1$ be dense in $X_0$ and $X_1$. If $b\in SV(0,\infty)$ satisfies
\begin{equation}\label{DC2}
\left\lVert t^{1-\theta-{1}/{q}} b(t)\min\left\lbrace 1,t^{-1}\right\rbrace\right\rVert_{q,(0,\infty)}<\infty,
\end{equation}
then 
\begin{equation}\label{embdg2}
(X_0',X_1')_{1-\theta,q\,'\!,b^{-1};J}\hookrightarrow (X_1,X_0)'_{\theta,q,b;K}.
\end{equation}
\end{theorem}
\bpr
By Theorem \ref{305}, condition \eqref{DC2} guarantees that the space $(X_0',X_1')_{1-\theta,q\,'\!,b^{-1};J}$ is a Banach space 
(which is intermediate between spaces $X_0'$ and $X_1'$). 

Since 
$$
t\min\{1,t^{-1}\}=\min\{1,t\} \quad \text{for all } \ t>0,
$$
it is clear that condition \eqref{DC2} is equivalent to
$$
\left\lVert t^{-\theta-{1}/{q}} b(t)\min\left\lbrace 1,t\right\rbrace\right\rVert_{q,(0,\infty)}<\infty.
$$
Consequently, by Theorem \ref{302}, the space $(X_1,X_0)_{\theta,q,b;K}$ is also a Banach, space 
(which is intermediate between spaces $X_0$ and $X_1$). 

Let $f'\in (X_0',X_1')_{1-\theta,q\,'\!,b^{-1};J}$. By Lemma \ref{JD}, there are $u_m'\in X_0'\cap X_1'$ for all $m\in\mathbb{Z}$ 
such that
\begin{equation}\label{BI}
f'=\sum_{m\in\mathbb{Z}}u_m'\hspace{1cm}(\text{convergence in }X_0'+X_1' ),
\end{equation}
\begin{equation}\label{BIE2N}
\lVert \{J(u_m',2^m,X_0',X_1')\}_{m\in\mathbb{Z}}\rVert_{\lambda_{1-\theta,q'\!,b^{-1}}}<\infty,
\end{equation} 
and 
\begin{equation}\label{BIE3N}
\|f'\|_{(X_0',X_1')_{1-\theta,q'\!,b^{-1};J}} \approx 
\inf \,\lVert \{J(u'_m,2^m)\}_{m\in \mathbb{Z}}\rVert_{\lambda_{1-\theta,q'\!,b^{-1}}},
\end{equation}
where the infimum extends over all sequences $\{u'_m\}_{m\in\mathbb{Z}}$ satisfying  \eqref{BI} and \eqref{BIE2N}.

By \eqref{2}, 
for any $f\in (X_1,X_0)_{\theta, q,b;K},$ 
\begin{align}\label{E2.1.2}
|\langle f',f\rangle|&\leq \sum_{m\in\mathbb{Z}} J(u_m',2^{-m};X_1',X_0')K(f,2^m;X_1,X_0)  \\
&= \sum_{m\in\mathbb{Z}} 2^{-m}J(u_m',2^{m};X_0',X_1')K(f,2^m;X_1,X_0).\nonumber
\end{align}
Since, cf. Remark \ref{dual of lamba space}, $\lambda_{\theta,q,b}$ and $\lambda_{1-\theta,q\,'\!,b^{-1}}$ are dual spaces via the duality $\sum\limits_{m\in\mathbb{Z}}2^{-m}\alpha_m\beta_m$, we get
$$
{\rm RHS \eqref{E2.1.2}} \leq \lVert \{J(u_m',2^{m};X_0',X_1')\}_{m\in\mathbb{Z}}\rVert_{\lambda_{1-\theta,q'\!,b^{-1}}}\,
\lVert \{K(f,2^m;X_1,X_0)\}_{m\in\mathbb{Z}}\rVert_{\lambda_{\theta,q,b}}.
$$
Now, by Lemma \ref{KD},
$$
\lVert \{K(f,2^m;X_1,X_0)\}_{m\in\mathbb{Z}}\rVert_{\lambda_{\theta,q,b}}\approx\lVert f\rVert_{(X_1,X_0)_{\theta,q,b;K}}
\quad \text{for all }\ f\in (X_1,X_0)_{\theta, q,b;K}.
$$
Together with \eqref{E2.1.2}, this implies that %
$$
\lVert f'\rVert_{(X_1,X_0)'_{\theta,q,b;K}}\lesssim 
 \lVert \{J(u_m',2^{m};X_0',X_1')\}_{m\in\mathbb{Z}}\rVert_{\lambda_{1-\theta,q'\!,b^{-1}}} 
\quad \text{for all }\ f'\in (X_0',X_1')_{1-\theta,q\,'\!,b^{-1};J}.
$$
Thus, using \eqref{BIE3N} and taking the infimum over all sequences $\{u_m\}_{m\in\mathbb{Z}}$ satisfying  \eqref{BI} and \eqref{BIE2N}, 
we arrive at
$$
\lVert f'\rVert_{(X_1,X_0)'_{\theta,q,b;K}}\lesssim \|f'\|_{(X_0',X_1')_{1-\theta,q'\!,b^{-1};J}}
\quad \text{for all }\ f'\in (X_0',X_1')_{1-\theta,q\,'\!,b^{-1};J}, 
$$
which means that embedding \eqref{embdg2} holds.
\epr

\smallskip
\begin{lemma}{\rm (\cite [Lemma 11.2]{OG})}\label{102} Assume that $1<q<\infty$.

\noindent
\!{\bf \,(i)} Let $b \in SV(0,\infty)$ be such that
\begin{equation}\label{prop_slow_var_funct_b1r}
\int_x^\infty t^{-1}b^{\,q}(t)\, dt <\infty  \quad\text{for all \ }x>0 \qquad \text{and\ \ }
\int_0^\infty t^{-1}b^{\,q}(t)\, dt =\infty.
\end{equation}
If 
\begin{equation}\label{def_slow_var_funct_a1r}
a(x):=b^{-q/q\,'}(x)\int_x^\infty t^{-1}b^{\,q}(t)\, dt  \quad\text{for all \ }x>0,
\end{equation}
then $a \in SV(0,\infty)$, 
\begin{equation}\label{prop_slow_var_funct_ar}
\int_0^x t^{-1}a^{-q\,'}(t)\, dt <\infty  \quad\text{for all \ }x>0 \qquad \text{and\ \ }
\int_0^\infty t^{-1}a^{-q\,'}(t)\, dt =\infty.
\end{equation}
Moreover,
\begin{equation}\label{def_slow_var_funct_br}
b(x)\approx a^{-q\,'\!/q}(x)\Big(\int_0^x t^{-1}a^{-q\,'}(t)\, dt \Big)^{-1} \quad\text{for a.a. \ }x>0
\end{equation}
and
\begin{equation}\label{103}
\Big(\int_0^x t^{-1}a^{-q\,'}(t)\, dt \Big)^{1/q\,'} \Big(\int_x^\infty t^{-1}b^{\,q}(t)\, dt \Big)^{1/q}
=\Big(\frac{1}{q\,'-1}\Big)^{1/q\,'} \quad\text{for all \ }x>0.
\end{equation}
\!{\bf \,(ii)} Let $a \in SV(0,\infty)$ be such that \eqref{prop_slow_var_funct_ar} holds. If 
\begin{equation}\label{def_slow_var_funct_br=}
b(x):=a^{-q\,'\!/q}(x)\Big(\int_0^x t^{-1}a^{-q\,'}(t)\, dt \Big)^{-1} \quad\text{for all \ }x>0,
\end{equation}
then $b \in SV(0,\infty)$ and \eqref{prop_slow_var_funct_b1r} is satisfied. Moreover, 
\begin{equation}\label{def_slow_var_funct_a1r=}
a(x)\approx b^{-q/q\,'}(x)\int_x^\infty t^{-1}b^{\,q}(t)\, dt  \quad\text{for a.a. \ }x>0,
\end{equation}
and 
\begin{equation}\label{103*}
\Big(\int_0^x t^{-1}a^{-q\,'}(t)\, dt \Big)^{1/q\,'} \Big(\int_x^\infty t^{-1}b^{\,q}(t)\, dt \Big)^{1/q}
=\Big(\frac{1}{q-1}\Big)^{1/q} \quad\text{for all \ }x>0.
\end{equation}
\end{lemma}

\smallskip
\begin{lemma}\label{hustota}
If $(X_0, X_1)$ is a compatible couple and $X_0\cap X_1$ is dense in $X_0$ and $X_1$, then $X_0\cap X_1$ is also dense in 
$X_0+X_1$.
\end{lemma}

\bpr
Let $f \in X_0+X_1$ and $f=f_0+f_1$, where $f_i \in X_i, i=1,2$, and let $\varepsilon >0$. Since 
$X_0\cap X_1$ is dense in $X_0$ and $X_1$, there are $g_i \in X_0\cap X_1$ satisfying 
$$
\|f_i - g_i\|< \varepsilon/2, \ \ \ i=1,2.
$$
Putting $g:=g_0+g_1$, we see that $g \in X_0\cap X_1,$ and 
$$
\|f-g\|_{X_0+X_1} = \inf \{\|h_0\|_{X_0}+\|h_1\|_{X_1}: f-g=h_0+h_1\},
$$
where the infimum extends over all representation $f-g=h_0+h_1$ of $f-g$ with $h_i \in X_i,$  $i=1,2.$ 

Since $f-g=(f_0-g_0) +(f_1-g_1)$, $f_i-g_i \in X_i,$ i=1,2, we get
$$
\|f-g\|_{X_0+X_1} \le \|f_0-g_0\|_{X_0} + \|f_1-g_1\|_{X_1}<\varepsilon/2+\varepsilon/2=\varepsilon.
$$
Consequently,  $X_0\cap X_1$ is dense in $X_0+X_1$.
\epr

\medskip
\section{Proof of Theorem \ref{DT0S}}\label{Pr1Main}

We claim that condition \eqref{DC1} of Lemma \ref{JDchara} 
is satisfied with $\theta=0$.  Indeed, 
\begin{align*}
\lVert t^{-1/q\,'}a^{-1}(t)\min\{1,t^{-1}\}\rVert_{q\,'\!,(0,\infty)}
&\leq\lVert t^{-{1}/{q\,'}}a^{-1}(t)\rVert_{q\,'\!,(0,1)}+ \lVert t^{-1-{1}/{q\,'}}a^{-1}(t)\rVert_{q\,'\!,(1,\infty)}\\
&=: I_1+I_2. 
\end{align*}
Note that $I_2<\infty$ by Lemma \ref{l2.6+} (i) and (iii). If $1<q<\infty$, then $I_1<\infty$ by 
Lemma~\ref{102}~(i). If $q=1$, then 
\begin{equation}\label{q=1}
a(x)=\int_x^\infty t^{-1}b(t)\,dt \quad \text { for all } \ x>0. 
\end{equation}
Thus,  
$$
I_1=\|a^{-1}(t)\|_{\infty, (0, 1)}=\Big( \int_1^\infty t^{-1}b(t)\,dt\Big)^{-1} <\infty,
$$
and our claim follows.

Applying Theorem \ref{E1} with $\theta=0,$ we obtain
\begin{equation}\label{0embdg1}
(X_0,X_1)'_{0,q,a;J}\hookrightarrow (X_1',X_0')_{1,q\,'\!,a^{-1};K}.
\end{equation}

Now we claim that also condition \eqref{DC2} of Theorem \ref{E2} with $\theta=0$ is satisfied. Indeed, 
\begin{equation}\label{nerov}
\lVert t^{1-1/q}b(t)\min\{1,t^{-1}\}\rVert_{q, (0,\infty)}
\leq\lVert t^{1-1/q}b(t)\rVert_{q,(0,1)}+ \lVert t^{-1/q}b(t)\rVert_{q,(1,\infty)}.
\end{equation}
The first term on {\rm RHS}\eqref{nerov} is finite by Lemma \ref{l2.6+} (iii), while the second one by \eqref{DT0A}.

Using Theorem \ref{E2} with $\theta=0$, we arrive at 
\begin{equation}\label{embdg23}
(X_1',X_0')_{1,q\,'\!,b^{-1};J}\hookrightarrow (X_0,X_1)'_{0,q,b;K}.
\end{equation}
Further, by Theorem \ref{equivalence theorem1}, 
\begin{equation}\label{L4.3}
{\rm RHS}\eqref{embdg23}=(X_0,X_1)'_{0,q,b;K}=(X_0,X_1)'_{0,q,a;J}={\rm LHS}\eqref{0embdg1}.
\end{equation}
Now we are going to show that 
{\rm RHS}\eqref{0embdg1}={\rm LHS}\eqref{embdg23}, 
i.e.,
\begin{equation}\label{LLL}
 (X_1',X_0')_{1,q\,'\!,a^{-1};K}=(X_1',X_0')_{1,q\,'\!,b^{-1};J}.                           
\end{equation}

To verify \eqref{LLL}, we consider two cases:

Let $1<q<\infty$. Then we apply Theorem \ref{ET1}  with $B:= a^{-1}, A:=b^{-1}$,  with $q$ replaced by $q\,'$, 
and with $(X_0, X_1)$ replaced by $(X'_1, X'_0)$.  
Note that then \eqref{prop_slow_var_funct_B} reads as \eqref{prop_slow_var_funct_ar}, 
and \eqref{def_slow_var_funct_A} reads as\eqref{def_slow_var_funct_br}. But \eqref{prop_slow_var_funct_ar} and 
\eqref{def_slow_var_funct_br} hold by Lemma \ref{102} (i). 
Thus, \eqref{LLL} is true by Theorem  \ref{ET1}.

If $q=1$, then $q\,'=\infty$ and the function $a$ is given by \eqref{q=1}. Now we make use of Theorem~\ref{ET1.1}
with $A:=b^{-1}, B:= a^{-1}$, with $q=\infty$ and with $(X_0, X_1)$ replaced by $(X'_1, X'_0)$. Note that then 
\eqref{c1.1} holds by our assumption \eqref{DT0A}, while 
\eqref{a.1} follows from \eqref{q=1}. Thus, \eqref{LLL} holds by Theorem~\ref{ET1.1}.

Using \eqref{L4.3}, \eqref{0embdg1}, \eqref{LLL}, and \eqref{embdg23}, we obtain 
\begin{align*}
(X_0,X_1)'_{0,q,b;K}&=(X_0,X_1)'_{0,q,a;J}\hookrightarrow (X_1',X_0')_{1,q\,'\!,a^{-1};K}\\
&=(X_1',X_0')_{1,q\,'\!,b^{-1};J} \hookrightarrow (X_0,X_1)'_{0,q,b;K}. 
\end{align*}
Consequently,
\begin{equation}\label{D}
(X_0,X_1)'_{0,q,b;K}=(X_1',X_0')_{1,q\,'\!,b^{-1};J}=(X_1',X_0')_{1,q\,'\!,a^{-1};K}.
\end{equation}
Furthermore, by Lemmas \ref{10} and \ref{01},
$$
(X_1',X_0')_{1,q\,'\!,b^{-1};J}=(X_0',X_1')_{0,q\,'\!,\tilde{b};J}\quad \text{and} \quad
(X_1',X_0')_{1,q\,'\!,a^{-1};K}=(X_0',X_1')_{0,q\,'\!,\tilde{a};K}. 
$$
Together with \eqref{D}, this gives \eqref{D2E}.
\hskip 8,6cm$\square$

\medskip
\section{Proof of Theorem \ref{DTJ1}}\label{Pr3Main}

Since $1<q<\infty$ and $a\in SV(0, \infty)$ satisfies 
\begin{equation}\label{prop_slow_var_funct_aJznovu}
\int_0^x t^{-1}a^{-q\,'}(t)\, dt <\infty  \quad\text{for all \ }x>0, \qquad 
\int_0^\infty t^{-1}a^{-q\,'}(t)\, dt =\infty,
\end{equation}
we can apply Theorem \ref{equivalence theorem} to get  
\begin{equation}\label{K_space'=J_space'}
(X_0, X_1)'_{0, q, a; J}=(X_0, X_1)'_{0, q, b; K}
\end{equation}
where
\begin{equation}\label{def_slow_var_funct_bznovu}
b(x) = a^{-q\,'\!/q}(x)\Big(\int_0^x t^{-1}a^{-q\,'}(t)\, dt \Big)^{-1} \quad\text{for all\ }x>0.
\end{equation}
Moreover, by Lemma \ref{102} (ii), 
\begin{equation}\label{prop_slow_var_funct_b1rznovu}
\int_x^\infty t^{-1}b^{\,q}(t)\, dt <\infty  \quad\text{for all \ }x>0, \qquad 
\int_0^\infty t^{-1}b^{\,q}(t)\, dt =\infty,
\end{equation}
and
$$
a(x)\approx b^{-q/q\,'}(x)\int_x^\infty t^{-1}b^{\,q}(t)\, dt,  \quad\text{for a.a. \ }x>0.
$$
Thus, by Theorem \ref{DT0S},
\begin{equation}\label{D2Eznovu}
(X_0,X_1)'_{0,q,b;K}=(X_0',X_1')_{0,q\,'\!,\tilde{b};J}=(X_0',X_1')_{0,\,q\,'\!,\tilde{a};K}, 
\end{equation} 
where where $\tilde{a}(x):=\frac{1}{a({1}/{x})}$ and $\tilde{b}(x):=\frac{1}{b({1}/{x})}$ for all $x>0$.
Now \eqref{K_space'=J_space'} and \eqref{D2Eznovu} imply that 
$$
(X_0, X_1)'_{0, q, a; J}=(X'_0, X'_1)_{0, q\,'\!,\tilde{a}; K}=(X'_0, X'_1)_{0, q\,'\!,\tilde{b}; J}.
$$
\hskip 14,9cm$\square$

\medskip
\section{Proof of Theorem \ref{DTJ11}}\label{Pr4Main}

Let the function $a$ satisfy the assumption of Theorem \ref{DTJ11} and let 
 $b \in  SV(0,\infty)$ be given by \eqref{def_slow_var_funct_b*1}.  Then, by Theorem \ref{equivalence theorem1**},
\begin{equation}\label{K_space=J_space1***}
(X_0, X_1)_{0, 1, a; J}=(X_0, X_1)_{0, 1, b; K}.
\end{equation}
 
Solving differential equation \eqref{def_slow_var_funct_b*1} and using the condition  
$a(\infty)=0,$  we arrive at 
\begin{equation}\label{aint}
a(x):=\int_x^{\infty}t^{-1}b(t)\,dt\quad \text{for all }x>0.
\end{equation}
Moreover, condition \eqref{prop_slow_var_funct_a*1} implies that the function $b$ satisfies condition \eqref{DT0A} with $q=1$ of 
Theorem \ref{DT0S}. Thus, on using this theorem, we arrive at 
\begin{equation}\label{D2E*}
(X_0,X_1)'_{0,1,b;K}=(X_0',X_1')_{0,\infty,\tilde{b};J}=(X_0',X_1')_{0,\infty,\tilde{a};K}, 
\end{equation} 
where $\tilde{b}(x):=\frac{1}{b({1}/{x})}$ and $\tilde{a}(x):=\frac{1}{a({1}/{x})}$ for all $x>0$.
Now, by \eqref{K_space=J_space1***} and \eqref{D2E*},
$$
(X_0, X_1)'_{0, 1, a; J}=(X_0',X_1')_{0,\infty,\tilde{a};K}=(X_0',X_1')_{0,\infty,\tilde{b};J},
$$
and the proof is complete.$\hskip 10,6cm\square$

\smallskip
\begin{remark}\label{veta 3.4}
{\rm In the proof of Theorem \ref{DTJ11} one can apply Theorem \ref{equivalence theorem1} instead of Theorem~\ref{equivalence theorem1**}. Indeed, making use of \eqref{aint} and the fact that condition \eqref{prop_slow_var_funct_a*1} implies that the function $b$ satisfies condition \eqref{prop_slow_var_funct_b1} with $q=1$ of Theorem \ref{equivalence theorem1}, one can get \eqref{K_space=J_space1***} from Theorem \ref{equivalence theorem1}.}
\end{remark}

\medskip
\section{Proofs of Theorem \ref{DT0S1} and Remark \ref{remark100*}}\label{Pr2Main} 

{\bf I.} Proof of Theorem \ref{DT0S1}. 
By Theorem \ref{KS},
\begin{equation}\label{*D2E1*}
(X_0,X_1)_{0,q,b;K}=(X_0, X_0+X_1)_{0,q,B;K}.
\end{equation}
The function $B$ satisfies (cf. \eqref{3131}, \eqref{DT0A1}, and \eqref{3141})
\begin{equation}\label{prop_slow_var_funct_B*}
\int_x^\infty t^{-1}B^{\,q}(t)\, dt <\infty  \quad\text{for all \ }x>0, \qquad 
\int_0^\infty t^{-1}B^{\,q}(t)\, dt =\infty.
\end{equation}
Thus, using \eqref{*D2E1*}, Lemma \ref{hustota}, Theorem \ref{DT0S} 
(with $(X_0, X_1)$ replaced by $(X_0,X_0+ X_1)$ and with $B, A$ instead of $b, a$), and \eqref{99}, we arrive at
\begin{align}\label{222}
(X_0,X_1)'_{0,q,b;K}&=(X_0, X_0+X_1)'_{0,q,B;K}\\
&=(X'_0, (X_0+X_1)')_{0,q\,'\!\!,\tilde{B};J}=(X'_0,(X_0+X_1)')_{0,q\,'\!\!,\tilde{A};K}\notag\\
&=(X'_0, X'_0\cap X_1')_{0,q\,'\!\!,\tilde{B};J}=(X'_0,X'_0\cap X_1')_{0,q\,'\!\!,\tilde{A};K}.\notag
\end{align}
Moreover, by Corollary  \ref{cor2/2},
\begin{equation}\label{*D2E1*c}
(X_0, X_1)_{0, q, b; K}=X_0 + (X_0,X_1)_{0, q, B; K}=X_0 + (X_0,X_1)_{0, q, A; J}.
\end{equation}
By Theorem \ref{302} (i),
$$ 
X_0 \cap X_1 \hookrightarrow (X_0,X_1)_{0,q,B;K}.
$$
Thus,
\begin{equation}\label{7} 
X_0 \cap X_1 \hookrightarrow X_0 \cap (X_0,X_1)_{0,q,B;K}.
\end{equation}
By our assumption, $X_0 \cap X_1$ is dense in $X_0$, and, by Theorem \ref{density theorem1} (with $b$ replaced by $B$), 
$X_0 \cap X_1$ 
is dense in $(X_0,X_1)_{0,q,B;K}.$ Together with embedding \eqref{7}, this shows that 
$X_0 \cap (X_0,X_1)_{0,q,B;K}$ is dense in $X_0$ and in $(X_0,X_1)_{0,q,B;K}$. Therefore, we can apply 
\eqref{99} (with $X_0, X_1$ replaced by $X_0, (X_0,X_1)_{0,q,B;K}$) to get that 
\begin{equation}\label{77}
(X_0 + (X_0,X_1)_{0, q, B; K})'=X'_0 \cap (X_0,X_1)'_{0, q, B; K}.
\end{equation}
Moreover, by Theorem \ref{DT0S} (with $B$ and $A$ instead of $b$ and $a$),
\begin{equation}\label{777}
(X_0,X_1)'_{0, q, B; K}=(X'_0,X'_1)_{0, q\,'\!\!,\tilde{B};J}.
\end{equation}
Using \eqref{*D2E1*c}, \eqref{77}, and  \eqref{777}, we arrive at 
\begin{equation}\label{7777}
(X_0, X_1)'_{0, q, b; K}= X'_0 \cap (X'_0,X'_1)_{0, q\,'\!\!,\tilde{B};J}.
\end{equation}

Let $1<q<\infty$. Then the function $A$ satisfies 
\begin{equation}\label{7A}
\int_0^x t^{-1}A^{-q\,'}(t)\, dt <\infty  \quad\text{for all \ }x>0, \qquad 
\int_0^\infty t^{-1}A^{-q\,'}(t)\, dt =\infty, 
\end{equation}
which follows from Lemma \ref{102} (i) (with $b$ and $a$ replaced by $B$ and $A$). 
By Theorem~\ref{305}~(i),
$$ 
X_0 \cap X_1 \hookrightarrow (X_0,X_1)_{0,q,A;J}.
$$
Thus,
\begin{equation}\label{71} 
X_0 \cap X_1 \hookrightarrow X_0 \cap (X_0,X_1)_{0,q,A;J}.
\end{equation}
By our assumption, $X_0 \cap X_1$ is dense in $X_0$, and, by Theorem \ref{density theorem} (with $A$ instead of~$a$), 
$X_0 \cap X_1$ is dense in $(X_0,X_1)_{0,q,A;J}.$ Together with embedding \eqref{71}, this shows that 
$X_0 \cap (X_0,X_1)_{0,q,A;J}$ is dense both in $X_0$ and in $(X_0,X_1)_{0,q,A;J}$. Therefore, we can apply 
\eqref{99} (with $X_0, X_1$ replaced by $X_0, (X_0,X_1)_{0,q,A;J}$) to get that 
\begin{equation}\label{771}
(X_0 + (X_0,X_1)_{0, q, A; J})'=X'_0 \cap (X_0,X_1)'_{0, q, A; J}.
\end{equation}
Moreover, by Theorem \ref{DTJ1} (with $A$ and $B$ instead of $a$ and $b$),
\begin{equation}\label{7772}
(X_0,X_1)'_{0, q, A; J}=(X'_0,X'_1)_{0, q\,'\!\!,\tilde{A};K}.
\end{equation}
Using \eqref{*D2E1*c}, \eqref{771}, and  \eqref{7772}, we arrive at 
\begin{equation}\label{77771}
(X_0, X_1)'_{0, q, b; K}= X'_0 \cap (X'_0,X'_1)_{0, q\,'\!\!,\tilde{A};K}.
\end{equation}

Let $q=1$. Then $q\,'=\infty$ and 
$$
A(x)=\int_x^\infty t^{-1} B(t)\, dt \quad \text {for all \ }  x>0.
$$
Hence, $A \in SV(0,\infty)\cap AC(0,\infty)$ satisfies
$$
A \text{\ is strictly decreasing},\quad A(0)=\infty, \quad A(\infty)=0,
$$
and 
$$
B(x):= -x \,A\, '(x) \quad\text{for a.a. }x>0\,.
$$
Consequently, by Theorem \ref{density theorem*} (with $A$ instead of $a$),
$X_0\cap X_1$ is dense in $(X_0, X_1)_{0, 1, A; J}$.
Therefore, as in the case that $1<q<\infty$, we obtain that \eqref{771} holds with $q=1$. 
Moreover, Theorem \ref{DTJ11} (with $A$, and $B$ instead of $a$, and $b$) implies that \eqref{7772} holds with $q=1$.
Thus, on using \eqref{*D2E1*c}, \eqref{771}, and  \eqref{7772}, we see that 
\eqref{77771} remains true if $q=1$.

Applying Theorem \ref{JS11} (with $X_0, X_1$ replaced by $X'_0, X'_1$ and 
with $q, a, \alpha,$ and $A$ replaced by $q\,', \tilde{b}, \tilde{\beta},$ and $ \tilde{B},$ 
respectively\ \footnotemark), 
\footnotetext{\ \ Note that \eqref{B1=} implies 
$$\tilde{A}(x):=(\tilde{B})^{-q/{q'}}(x) \left(\int_0^x t^{-1}(\tilde{B})^{-q}(t)\,dt\right)^{-1} \quad\text{for all \ }x>0,$$
which is \eqref{Bfce111} with the above mentioned replacement in Theorem \ref{JS11}.}
we obtain 
$$
(X'_0, X'_1)_{0,q\,'\!\!,\tilde{b};J}=(X'_0, X'_0\cap X'_1)_{0,q\,'\!\!,\tilde{B};J}=(X'_0, X'_0\cap X'_1)_{0,q\,'\!\!,\tilde{A};K},
$$
which, together with \eqref{222}, \eqref{7777} and \eqref{77771}, implies that  \eqref{D2E1} holds.
\hskip 3,37cm$\square$

\smallskip
{\bf II.} Proof of Remark \ref{remark100*}. Equality \eqref{101*} holds by \eqref{101}, while \eqref{1002*}-\eqref{1004*} follow from Remark \ref{remark100J} (used with $X'_0, X'_1$ instead of $X_0, X_1$ and 
with $q\,'\!\!, \tilde{b}, \tilde{B},$ and $ \tilde{A}$ instead of 
$q, a, A,$ and $B$, respectively).
\hskip 11,28cm$\square$

\smallskip

\section{Proofs of Theorem \ref{DTJ111} and Remark \ref{remark100J*}} \label{Pr5Main}

{\bf I.} Proof of Theorem \ref{DTJ111}. 
By Theorem \ref{JS11},
\begin{equation}\label{corres2*}
(X_0, X_1)_{0, q, a; J}=(X_0, X_0\cap X_1)_{0,q,A;J}.
\end{equation}
The function $A$ satisfies (cf. \eqref{Afce}, \eqref{akon}, and \eqref{alfa})
\begin{equation}\label{A111}
\int_0^x t^{-1}A^{-q\,'}(t)\, dt <\infty\quad\text{for all \ }x>0, \qquad 
\int_0^\infty t^{-1}A^{-q\,'}(t)\, dt =\infty.
\end{equation}
Thus, using \eqref{corres2*}, Theorem \ref{DTJ1} (with $(X_0,X_1)$ replaced by $(X_0,X_0 \cap X_1)$ and 
with $A, B$ instead of $a, b$), and \eqref{99}, 
we arrive at
\begin{align}\label{*corres2*}
(X_0, X_1)'_{0, q, a; J}&=(X_0, X_0\cap X_1)'_{0,q,A;J}\\
&=(X'_0, (X_0\cap X_1)')_{0,q\,'\!\!,\tilde{A};K}=(X'_0, (X_0\cap X_1)')_{0, q\,'\!\!, \tilde{B}; J}\notag\\
&=(X'_0, X'_0+X'_1)_{0,q\,'\!\!,\tilde{A};K}=(X'_0, X'_0+ X'_1)_{0, q\,'\!\!, \tilde{B}; J}\,.\notag
\end{align}
Moreover, by Corollary \ref{cor2},
\begin{equation}\label{*corres2}
(X_0, X_1)_{0, q, a; J}=X_0 \cap (X_0,X_1)_{0, q, A; J}=X_0 \cap (X_0,X_1)_{0, q, B; K}.
\end{equation}
As in the proof of Theorem \ref{DT0S1}, one can show that $X_0 \cap (X_0,X_1)_{0,q,A;J}$ is dense both in $X_0$ and 
in $(X_0,X_1)_{0,q,A;J}$. Therefore, applying  
\eqref{99} (with $X_0, X_1$ replaced by $X_0, (X_0,X_1)_{0,q,A;J}$), we arrive at  
\begin{equation}\label{44}
(X_0 \cap (X_0,X_1)_{0, q, A; J})'=X'_0 + (X_0,X_1)'_{0, q, A; J}.
\end{equation}
Moreover, by Theorem \ref{DTJ1} (with $A$ and $B$ instead of $a$ and $b$), \eqref{7772} holds. 
Using \eqref{*corres2}, \eqref{44}, and  \eqref{7772}, we arrive at 
\begin{equation}\label{777714}
(X_0, X_1)'_{0, q, a; J}= X'_0 + (X'_0,X'_1)_{0, q\,'\!\!,\tilde{A};K}.
\end{equation}
Making use of \eqref{A111}, \eqref{Bfce}, and Lemma \ref{102} (ii) (with $a$ and $b$ replaced by $A$ and $B$), 
we obtain that \eqref{prop_slow_var_funct_B*} holds. As in the proof of Theorem \ref{DT0S1}, one can show that 
$X_0 \cap (X_0,X_1)_{0,q,B;K}$ is dense both in $X_0$ and in $(X_0,X_1)_{0,q,B;K}$. Therefore, applying  
\eqref{99} (with $X_0, X_1$ replaced by $X_0, (X_0,X_1)_{0,q,B;K}$), we get
\begin{equation}\label{774}
(X_0 \cap (X_0,X_1)_{0, q, B; K})'=X'_0 + (X_0,X_1)'_{0, q, B; K}.
\end{equation}
Moreover, by Theorem \ref{DT0S} (with $B$ and $A$ instead of $b$ and $a$), \eqref{777} holds. Relations 
\eqref{*corres2}, \eqref{774}, and \eqref{777} imply that 
\begin{equation}\label{7744}
(X_0, X_1)'_{0, q, a; J}=X'_0+(X'_0,X'_1)_{0, q\,'\!\!, \tilde{B}; J}\,.
\end{equation}
Applying Theorem \ref{KS} (with $X_0, X_1$ replaced by $X'_0, X'_1$ and 
with $q, b, \beta,$ and $B$ replaced by $q\,', \tilde{a}, \tilde{\alpha},$ and $ \tilde{A},$ 
respectively\ \footnotemark), 
\footnotetext{\ \ 
Note that \eqref{Bfce} implies  
$$\tilde{B}(x):=(\tilde{A})^{-{q'\!}/{q}}(x)\int_x^{\infty}t^{-1}(\tilde{A})^{q'}(t)\,dt\quad \text{for all }x>0,$$
which is \eqref{B1=*} with the above mentioned replacement in Theorem \ref{KS}.}
we obtain 
$$
(X'_0, X'_1)_{0,q\,'\!\!,\tilde{a};K}=(X'_0, X'_0+X'_1)_{0,q\,'\!\!,\tilde{A};K}=(X'_0, X'_0+X'_1)_{0,q\,'\!\!,\tilde{B};J},
$$
which, together with \eqref{*corres2*}, \eqref{777714}, and \eqref{7744} shows that  \eqref{dual} holds.
\hskip 3,77cm$\square$

\smallskip
{\bf II.} Proof of Remark \ref{remark100J*}. Equality \eqref{101J*} holds by \eqref{101J}, 
while \eqref{1002J*}-\eqref{1004J*} follow from Remark \ref{remark100} (used with $X'_0, X'_1$ instead of $X_0, X_1$ and 
with $q\,'\!\!, \tilde{a}, \tilde{A},$ and $ \tilde{B}$ instead of $q, b, B,$ and $A$, respectively).
\hskip 11,3cm$\square$
\medskip

\section{Proofs of Theorem \ref{DTJ111.1} and  and Remark \ref{remark100J*4}}\label{Pr6Main}

{\bf I.} Proof of Theorem \ref{DTJ111.1}. By Theorem \ref{equivalence theorem1**a},
\begin{equation}\label{corres2*.2}
(X_0, X_1)_{0, 1, a; J}=(X_0, X_0\cap X_1)_{0,1,A;J}.
\end{equation}
The function $A \in SV(0,\infty) \cap AC(0,\infty)$ satisfies 
\begin{equation}\label{A111.2}
A \text{  is strictly decreasing}, \qquad A(0)=\infty, \qquad A(\infty)=0.
\end{equation}
Thus, using \eqref{corres2*.2}, Theorem \ref{DTJ11} with $(X_0,X_1)$ replaced by $(X_0,X_0 \cap X_1)$ and 
with $A, B$ instead of $a, b$), and \eqref{99}, we arrive at 
\begin{align}\label{*corres2*2}
(X_0, X_1)'_{0, 1, a; J}&=(X_0, X_0\cap X_1)'_{0,1,A;J}\\
&=(X'_0, (X_0\cap X_1)')_{0, \infty,\tilde{A};K}=(X'_0, (X_0\cap X_1)')_{0, \infty, \tilde{B}; J}\notag\\
&=(X'_0, X'_0+X'_1)_{0, \infty,\tilde{A};K}=(X'_0, X'_0+ X'_1)_{0, \infty, \tilde{B}; J}\,.\notag
\end{align}
Moreover, by Corollary \ref{cor21},
\begin{equation}\label{*corres2.2}
(X_0, X_1)_{0, 1, a; J}=X_0 \cap (X_0,X_1)_{0, 1, A; J}=X_0 \cap (X_0,X_1)_{0, 1, B; K}.
\end{equation}
By Theorem~\ref{305}~(i),
$$ 
X_0 \cap X_1 \hookrightarrow (X_0,X_1)_{0,1,A;J}.
$$
Thus,
\begin{equation}\label{714} 
X_0 \cap X_1 \hookrightarrow X_0 \cap (X_0,X_1)_{0,1,A;J}.
\end{equation}
By our assumption, $X_0 \cap X_1$ is dense in $X_0$, and, by Theorem \ref{density theorem*} (with $A$ instead of~$a$), 
$X_0 \cap X_1$ is dense in $(X_0,X_1)_{0,q,A;J}.$ Together with embedding \eqref{714}, this shows that 
$X_0 \cap (X_0,X_1)_{0,1,A;J}$ is dense both in $X_0$ and in $(X_0,X_1)_{0,1,A;J}$. Therefore, we can apply 
\eqref{99} (with $X_0, X_1$ replaced by $X_0, (X_0,X_1)_{0,q,A;J}$) to get that 
\begin{equation}\label{7711}
(X_0 \cap (X_0,X_1)_{0, 1, A; J})'=X'_0 + (X_0,X_1)'_{0, 1, A; J}.
\end{equation}
Moreover, by Theorem \ref{DTJ11} (with $A$ and $B$ instead of $a$ and $b$),
\begin{equation}\label{7773}
(X_0,X_1)'_{0, 1, A; J}=(X'_0,X'_1)_{0, \infty,\tilde{A};K}.
\end{equation}
Using \eqref{*corres2.2}, \eqref{7711}, and  \eqref{7773}, we arrive at 
\begin{equation}\label{777711}
(X_0, X_1)'_{0, 1, A; J}= X'_0 + (X'_0,X'_1)_{0, \infty,\tilde{A};K}.
\end{equation}

Solving differential equation \eqref{Bfce.1} and using the condition $A(\infty)=0$, we obtain that
\begin{equation}\label{aint1}
A(x):=\int_x^{\infty}t^{-1}B(t)\,dt\quad \text{for all }x>0,
\end{equation}
which, together with the condition $A(0)=\infty$, implies that 
\begin{equation}\label{prop_slow_var_funct_B*1}
\int_x^\infty t^{-1}B(t)\, dt <\infty  \quad\text{for all \ }x>0, \qquad 
\int_0^\infty t^{-1}B(t)\, dt =\infty.
\end{equation}
By Theorem \ref{302} (i),
$$ 
X_0 \cap X_1 \hookrightarrow (X_0,X_1)_{0,1,B;K}.
$$
Thus,
\begin{equation}\label{7*} 
X_0 \cap X_1 \hookrightarrow X_0 \cap (X_0,X_1)_{0,1,B;K}.
\end{equation}
By our assumption, $X_0 \cap X_1$ is dense in $X_0$, and, by Theorem \ref{density theorem1} (with $B$ instead of $b$), 
$X_0 \cap X_1$ is dense in $(X_0,X_1)_{0,1,B;K}.$ Together with embedding \eqref{7*}, this shows that 
$X_0 \cap (X_0,X_1)_{0,1,B;K}$ is dense both in $X_0$ and in $(X_0,X_1)_{0,1,B;K}$. Therefore, we can apply 
\eqref{99} (with $X_0, X_1$ replaced by $X_0, (X_0,X_1)_{0,1,B;K}$) to get that 
\begin{equation}\label{778}
(X_0 \cap (X_0,X_1)_{0, 1, B; K})'=X'_0 + (X_0,X_1)'_{0, 1, B; K}.
\end{equation}
Moreover, by Theorem \ref{DT0S} (with $B$ and $A$ instead of $b$ and $a$),
\begin{equation}\label{7779}
(X_0,X_1)'_{0, 1, B; K}=(X'_0,X'_1)_{0, \infty,\tilde{B};J}.
\end{equation}
Using \eqref{*corres2.2}, \eqref{778}, and  \eqref{7779}, we arrive at 
\begin{equation}\label{777712}
(X_0, X_1)'_{0, 1, a; J}= X'_0 + (X'_0,X'_1)_{0, \infty,\tilde{B};J}.
\end{equation}
Applying Theorem \ref{*equivalence theorem1*} (with $X_0, X_1$ replaced by $X'_0, X'_1$ and 
with $b, \beta, B,$ and $A$ replaced by $\tilde{a}, \tilde{\alpha}, \tilde{A},$ and 
$\tilde{B}$, respectively\ \footnotemark),  
\footnotetext{\ \ 
Note that \eqref{Bfce.1} implies  
$$\tilde{B}(x):= \frac{(\tilde{A})^2(x)}{x (-(\tilde{A})'(x))} \quad\text{for a.a. }x>0,$$ 
which is \eqref{difeq} with the above mentioned replacement in Theorem \ref{*equivalence theorem1*}.}
$$
(X'_0, X'_1)_{0,\infty,\tilde{a};K}=(X'_0, X'_0+X'_1)_{0,\infty,\tilde{A};K}=(X'_0, X'_0+X'_1)_{0,\infty,\tilde{B};J},
$$
which, together with \eqref{*corres2*2}, \eqref{777711}, and \eqref{777712}, shows that  \eqref{vnoreni.1} holds.
\hskip 2,88cm$\square$

\smallskip

{\bf II.} Proof of Remark \ref{remark100J*4}. Equality \eqref{101J*4} holds by \eqref{101J1}, 
while \eqref{1002J*4}-\eqref{1004J*4} follow from Remark \ref{remark100*4} 
 (with $X_0, X_1$ replaced by $X'_0, X'_1$ and 
with $b, \beta, B,$ and $A$ replaced by $\tilde{a}, \tilde{\alpha}, \tilde{A},$ and 
$\tilde{B}$, respectively).
\hskip 10,01cm$\square$

\medskip
\section{Concluding remarks}\label{CR}

Note that in \cite{CP} discrete versions of $K$- and $J$-spaces are defined and then in \cite [Theorems 3.1 and 3.2]{CP} 
duals of these spaces are described. Thus, there is another way how to prove some results mentioned in Section \ref{sectionmain} 
of our paper. For example, to prove the first equality mentioned in \eqref{D2E} of Theorem \ref{DT0S}, one can proceed as follows:

\noindent 
1) Use Lemma \ref{KD} to rewrite the given space $(X_0,X_1)_{0,q,b;K}$ in the discrete form.

\noindent 
2) Apply results \cite [Theorems 3.1]{CP} to determine the dual of the given discrete $K$-space. 
(Note that the dual is expressed as a discrete $J$-space.)

\noindent 
3) Use Lemma \ref{JD} to rewrite this discrete $J$-space in the continuous form.

Note that such a method was used in \cite [Theorem 5.6]{CS}, where the particular case of Theorem \ref{DT0S}, with $b$ of a logarithmic form, was proved.
\medskip

\bibliographystyle{alpha}

\begin{thebibliography}{CDSFM15}

\bibitem[AEEK11]{AEEK11}
I.~Ahmed, D. E. Edmunds, W.D. Evans, and G. E. Karadzhov. 
\newblock Reiteration theorems for the K-interpolation method in limiting cases.
\newblock {\em Math. Nachr.}, 284(4):421--442, 2011.

\bibitem[AFH20]{AFH20}
I.~Ahmed, A. Fiorenza, and A. Hafeez. 
\newblock Some interpolation formulae for grand and small Lorentz spaces.
\newblock {\em Mediterr. J. Math.}, 17(2), Art. 57, 21 pp., 2020. 


\bibitem[BS88]{BS:IO}
C.~Bennett and R.~Sharpley.
\newblock {\em Interpolation of Operators}, volume 129 of {\em Pure and Applied
  Mathematics}.
\newblock Academic Press, New York, 1988.

\bibitem[BL76]{BL}
J.~Bergh and J.~L\" ofstr\" om.
\newblock {\em Interpolation spaces. An introduction}, volume 223 of {\em Grundlehren der Mathematischen Wissenschaften.}
\newblock Springer-Verlag, Berlin-New York, 1976.  

\bibitem[BCFC20]{BCFC20}
B.~F. Besoy, F.~Cobos, and L. M. Fernández-Cabrera.
\newblock On the Structure of a Limit Class of Logarithmic Interpolation Spaces.
\newblock {\em  Mediterr. J. Math.}, 17(5), Art. 168, 17 pp., 2020.  




\bibitem[BGT87]{BGT87:RV}
N.~H. Bingham, C.~M. Goldie, and J.~L. Teugels.
\newblock {\em Regular Variation}.
\newblock Cambridge University Press, Cambridge, 1987.

\bibitem[BK91]{BK}
Yu.~A. Brudnyi and N.~Ya.~Krugljak.
\newblock {\em Interpolation Functors and Interpolation Spaces, vol. 1}, volume 47 of {\em North-Holland Mathematical Library}.
\newblock North-Holland Publishing Co., Amsterdam, 1991. 

\bibitem[BB67]{BB}
P.~L. Butzer and H.~Berens.
\newblock {\em Semi-groups of operators and approximation}, volume 145 of  {\em Die Grundlehren der 
mathematischen Wissenschaften}.
\newblock Springer-Verlag New York Inc., New York, 1967. 












\bibitem[CD15]{COBOS201543}
F.~Cobos and {\'O}.~Dom\'{\i}nguez.
\newblock Approximation spaces, limiting interpolation and {B}esov spaces.
\newblock {\em J. Approx. Theory}, 189:43--66, 2015.


\bibitem[CD16]{CD3}
F.~Cobos and {\'O}.~Dom\'{i}nguez.
\newblock On the relationship between two kinds of {B}esov spaces with
  smoothness near zero and some other applications of limiting interpolation.
\newblock {\em J. Fourier Anal. Appl.}, 22(5):1174--1191, 2016.

\bibitem[CDT16]{CDT3}
F.~Cobos,  {\'O}.~Dom\'{i}nguez, and H. Triebel.
\newblock Characterizations of logarithmic Besov spaces in terms of differences, Fourier-analytical decompositions, 
wavelets and semi-groups.
\newblock {\em J. Funct. Anal.}, 270:4386--4425, 2016.



\bibitem[CFCKU09]{CFCKU}
F. Cobos, L. M. Fern\' andez-Cabrera, T. K\" {u}hn, and T. Ullrich.
\newblock On an extreme class of real interpolation spaces.
\newblock {\em J. Funct. Anal.} 256(7):2321--2366, 2009.

\bibitem[CFCG24]{CFCG24}
F. Cobos, L. M. Fern\' andez-Cabrera, and M. Grover.
\newblock Measure of non-compactness and limiting interpolation
with slowly varying functions.
\newblock {Banach J. Math. Anal.} 18:25, 2024.
https://doi.org/10.1007/s43037-024-00335-z

\bibitem[CFCM15]{CFCM15}
F. Cobos, L. M. Fern\' andez-Cabrera, and A.  Mart{\' \i}nez.
\newblock On a paper of Edmunds and Opic on limiting interpolation of compact 
operators between $L_p$ spaces.
\newblock {\em Math. Nachr.} 288(2--3):167--175, 2015.

\bibitem[CFCM10]{CFCM10}
F. Cobos, L. M. Fern\' andez-Cabrera, and M. Masty{\l}o.
\newblock Abstract limit $J$-spaces.
\newblock {\em J. Lond. Math. Soc.}  82(2):501--525, 2010. 

\bibitem[CFC19]{CFC19}
F. Cobos and L. M. Fern\' andez-Cabrera.
\newblock Duality for logarithmic interpolation spaces and applications to Besov spaces.
\newblock {\em Banach Center Publ.}  119:109–122. 
\newblock Polish Academy of Sciences, Institute of Mathematics, Warsaw, 2019. 

\bibitem[CK11]{CK}
F.~Cobos, T.~K\" {u}hn.
\newblock Equivalence of K- and J-methods for limiting real interpolation spaces.
\newblock {\em J. Funct. Anal.} 261(12):3696--3722, 2011. 

\bibitem[CS15]{CS}
F.~Cobos, A.~Segurado.
\newblock Description of logarithmic interpolation spaces by means of the J-functional and applications. 
\newblock {\em J. Funct. Anal.} 268(10):2906--2945, 2015.

\bibitem[CP81]{CP}
M.~Cwikel, J.~Peetre.
\newblock Abstract K and J spaces.
\newblock {\em J. Math. pures et appl.} 60:1--50, 1981.



\bibitem[Dok91]{Dok91}
 R.~Ya. Doktorskiĭ.
\newblock Reiterative relations of the real interpolation method {\rm (}Russian{\rm )}. 
\newblock  {\em Dokl. Akad. Nauk SSSR}  321(2):241--245, 1991. 

\bibitem[Dok18]{Dok18}
 R.~Ya. Doktorskiĭ.
\newblock An application of limiting interpolation to Fourier series theory.
{\em The diversity and beauty of applied operator theory.}  179--191, volume 268 of {\em  Oper. Theory Adv. Appl.}
 \newblock Birkhäuser/Springer, Cham, 2018.




\bibitem[DT23]{ODST2019MAMS}
{\'O}.~Dom\'{i}nguez and S.~Tikhonov.
\newblock Function spaces of logarithmic smoothness: embeddings and characterizations.
\newblock {\em Mem. Amer. Math. Soc.}, 282(1393), vii+166 pp., 2023.


\bibitem[EE04]{EEv04:HOFSE}
D.~E. Edmunds and W.~D. Evans.
\newblock {\em Hardy {O}perators, {F}unction {S}paces and {E}mbeddings}.
\newblock Springer-Verlag, Berlin, Heidelberg, 2004.


\bibitem[EKP00]{EKP:OISRIQ}
D.~E. Edmunds, R.~Kerman, and L.~Pick.
\newblock Optimal {S}obolev {I}mbeddings {I}nvolving
  {R}earrangement-{I}nvariant {Q}uasinorms.
\newblock {\em J. Funct. Anal.}, 170:307--355, 2000.

\bibitem[EO14]{EO14}
D.~E. Edmunds and B. Opic.
\newblock Limiting variants of Krasnosel'ski\u{\i}'s compact interpolation theorem.
\newblock {\em J. Funct. Anal.}, 266(5): 3265--3285, 2014.

\bibitem[EO00]{EO00:RILFR}
W.~D. Evans and B.~Opic.
\newblock Real interpolation with logarithmic functors and reiteration.
\newblock {\em Canad. J. Math.}, 52:920--960, 2000.

\bibitem[EOP02]{EOP02:RILF}
W.~D. Evans, B.~Opic, and L.~Pick.
\newblock Real {I}nterpolation with {L}ogarithmic {F}unctors.
\newblock {\em J. Inequal. Appl.}, 7(2):187--269, 2002.

\bibitem[FMS12]{FMS12}
P. Fern\' andez-Mart{\'\i}nez and T. Signes.
\newblock Real interpolation with symmetric spaces and slowly varying functions.
\newblock {\em Q. J. Math.}, 63(1):133--164, 2012.

\bibitem[FMS14]{FMS14}
P. Fern\' andez-Mart{\'\i}nez and T. Signes.
\newblock Reiteration theorems with extreme values of parameters.
\newblock {\em  Ark. Mat.}, 52(2):227--256, 2014.

\bibitem[FMS15]{FMS15}
P. Fern\' andez-Mart{\'\i}nez and T. Signes.
\newblock Limit cases of reiteration theorems.
\newblock {\em Math. Nachr.}, 288(1):25--47, 2015.



\bibitem[FFGKR18]{FFGKR18}
A.~Fiorenza, M.~R. Formica, A.~Gogatishvili, T. Kopaliani, and J. M. Rakotoson.
\newblock Characterization of interpolation between grand, small or classical Lebesgue spaces.
\newblock {\em Nonlinear Anal.}, 177(B):422--453, 2018.



\bibitem[GNO10]{GNO10:PotAnal}
A.~Gogatishvili, J.~S. Neves, and B.~Opic.
\newblock Optimal embeddings of {B}essel-potential-type spaces into generalized
  {H}\"older spaces involving $k$-modulus of smoothness.
\newblock {\em Potential Anal.}, 32(3):201--228, 2010.

\bibitem[GOT05]{GOT2002Lrrinwsvf}
A.~Gogatishvili, B.~Opic, and W.~Trebels.
\newblock Limiting reiteration for real interpolation with slowly varying
  functions.
\newblock {\em Math. Nachr.}, 278:86--107, 2005.






\bibitem[GM86]{GM86}
M. E. Gomez, M. Milman.
\newblock Extrapolation spaces and almost-everywhere convergence of singular integrals. 
\newblock {\em J. London Math. Soc.} 34(2):305--316, 1986.
 


\bibitem[Gu78]{Gu78}
J. Gustavsson.
\newblock A function parameter in connection with interpolation of Banach spaces.
\newblock {\em Math. Scand.},  42(2):289--305, 1978.

\bibitem[KPS78]{KPS78}
S.~G. Kre{\u \i}n, Ju.~I. Petunin and E.~M. Sem\" {e}nov.
\newblock {\em Interpolation of linear operators {\rm (}Russian{\rm )}. } 
 {\em   Nauka.}
\newblock Moscow, 1978.

\bibitem[Li61]{Li61}
J.-L. Lions. 
\newblock  Sur les espaces d'interpolation: dualit\' e.
\newblock {\em Math. Scand.,} 9:147--177, 1961.


\bibitem[LP64]{LP64}
J.-L. Lions and J. Peetre. 
\newblock  Sur une classe d´espaces d´interpolation.
\newblock {\em Inst. Hautes ´Etudes Sci. Publ. Math.,} 19:5--68, 1964.


\bibitem[L18]{L18}
A. Lunardi.
\newblock {\em Interpolation theory. Third edition.}, volume 16 of  {\em Appunti.  
  Scuola Normale Superiore di Pisa {\rm(}Nuova Serie{\rm)} {\rm[}Lecture Notes. Scuola Normale Superiore di Pisa {\rm (}New Series{\rm )}{\rm]}}.
\newblock Edizioni della Normale, Pisa, 2018.





\bibitem[Mal84]{Mal84}
L.~Maligranda.
\newblock The K-functional for symmetric spaces, 
volume 1070 of
  {\em Lecture Notes in Mathematics}.
\newblock Springer-Verlag, Berlin, 1984: 169--182.

\bibitem[Mar00]{VojislavMaric:RVDE00}
V.~Mari\'c.
\newblock {\em Regular Variation and Differential Equations}, volume 1726 of
  {\em Lecture Notes in Mathematics}.
\newblock Springer-Verlag, Berlin, 2000.



\bibitem[Mil94]{Mil94}
M.~Milman.
\newblock {\em  Extrapolation and optimal decompositions with applications to analysis.}, 
volume~1580 of {\em  Lecture Notes in Mathematics.}
\newblock Springer-Verlag, Berlin, 1994.

\bibitem[N02]{Nev02:LKSBRPE}
J.~S. Neves.
\newblock Lorentz-{K}aramata spaces, {B}essel and {R}iesz potentials and
  embeddings.
\newblock {\em Dissertationes Math. {\rm (}Rozprawy Mat.{\rm )}}, 405:46 pp., 2002.

\bibitem[NO20]{NO20}
J.~S. Neves and B. Opic.
\newblock Optimal local embeddings of Besov spaces involving only slowly varying smoothness.
\newblock {\em J. Approx. Theory}, 254 (105393), 25 pp., 2020.

\bibitem[OG23]{OG} 
B. Opic and M. Grover. 
\newblock Description of $K$-spaces by means of $J$-spaces and the reverse problem in the limiting interpolation.
\newblock {\em Math. Nachr.},  296:4002–4031, 2023. 
 
\bibitem[OG24]{O24} 
B. Opic and M. Grover. 
\newblock Relationship between limiting $K$-spaces and $J$-spaces in the real interpolation,
\newblock {arXiv:2501.08965 [math.FA]}.

 \noindent	
https://doi.org/10.48550/arXiv.2501.08965

\bibitem[OK90]{OK90:HTI}
B.~Opic and A.~Kufner.
\newblock {\em Hardy-type inequalities}.
\newblock Pitman Research Notes in Math. Series 219, Longman Sci. \& Tech.,
  Harlow, 1990.
\bibitem[P86]{P86}

L.E. Persson.
\newblock Interpolation with a parameter function.
\newblock {\em Math. Scand.}, 59:199--222, 1986.






\bibitem[Tri78]{Tri78}
H.~Triebel.
\newblock {\em Interpolation theory, function spaces, differential operators.}, 
volume 18 of {\em  North-Holland Mathematical Library.}
\newblock North-Holland Publishing Co., Amsterdam-New York, 1978.

\bibitem[Zyg57]{Zyg57:TS}
A.~Zygmund.
\newblock {\em Trigonometric Series}, volume~I.
\newblock Cambridge University Press, Cambridge, 1957.

\end{thebibliography}

\def\cprime{$'$}

\end{document}